\documentclass[11 pt]{article}

 \pdfoutput=1
\usepackage{amsmath}
\usepackage{epsfig}
\usepackage{graphicx,color}
\usepackage{amsfonts}
\usepackage{amssymb}
\usepackage{verbatim}
\usepackage{float}
\usepackage{subfig}
\usepackage{tabularx}
\usepackage{array}
\usepackage{color}

\restylefloat{figure}

\numberwithin{equation}{section}

\newtheorem{defn}{Definition}[section]

\newtheorem{lem}[defn]{Lemma}

\makeatletter

\newcommand*{\Rom}[1]{\expandafter\@slowromancap\romannumeral #1@}
\makeatother

\begin{document}

{
\title{Influence of awareness that results from direct experience on the spread of epidemics}

\author{Ying~Xin \\ Department of Mathematics\\  Ohio University}

\date{\today}

\maketitle
}\noindent

\begin{abstract}

Here we study ODE epidemic models with spread of awareness, assuming that a certain proportion of the hosts will become aware of the ongoing outbreak upon recovery.  This study builds on W. Just and J. Salda\~{n}a's work in \cite{WJ}, and is conducted under the same framework, while addressing the influence of the awareness gained from direct experience of the disease.

In \cite{WJ}, the authors investigated the question whether preventive behavioral response triggered by awareness of the infection is sufficient to prevent future flare-ups from low endemic levels if awareness decays over time.  They showed that if all the hosts experienced infection return directly to the susceptible compartment upon recovery, such oscillations are ruled out in Susceptible-Aware-
Infectious-Susceptible models with a single compartment of aware hosts, but
can occur if two distinct compartments of aware hosts who differ
in their willingness to alert other susceptible hosts are considered.  Qualitatively, the models studied here produce the same results when we assume that recovery from the disease may or even will convey awareness from direct experience.

\end{abstract}

\section{Introduction}\label{s1}

Behavioral responses to an infectious disease are based on awareness that can result either from direct experience or information about an ongoing outbreak.

In \cite{WJ}, the authors built reactive SAIS and SAUIS models to study the influence of such behavioral responses arising from awareness that decays over time, on epidemic spreading.  They have mainly focused on the question under what circumstances
a behavioral response that is induced by awareness can be an effective control
measure.  That is, whether such models would predict a lowered epidemic threshold and
whether the response would prevent future flare-ups from low endemic levels.  A detailed discussion of the motivation for such models and the broader literature on this subject can be found in \cite{WJ}.  However, in these models, it is assumed that all the hosts who experienced infection will return to the susceptible compartment directly upon recovery.  That is, possible awareness that results from direct experience is totally ignored in these models, which is not very realistic.  

Here with the same questions in mind, the SAIAS and SAUIUAS models are defined and explored.  These models assume that a certain proportion of the hosts will become aware (or unwilling when applicable) upon recovery, instead of going back to the susceptible compartment directly.  We will show that these SAIAS and SAUIUAS models show the same qualitative behaviors as the SAIS and SAUIS models of \cite{WJ}.   Specifically, Section 2 shows that oscillations are ruled out in the SAIAS models regardless of the level of awareness gained from direct experience.  Section 3 shows that future flare-ups from low endemic levels are still possible in the SAUIUAS models.  However, within certain regions of the parameter space, originally possible oscillations can be eliminated by higher levels of awareness gained from direct experience. 

Most of our calculations  and the presentation of the material closely follow the more detailed exposition in \cite{WJ}.

\section{Reactive SAIAS models}\label{SAIS-Non}

\subsection{The model}

Similar to the SAIS models in \cite{WJ}, an SAIAS model has three compartments: S (susceptible), A (aware) and I (infectious). Susceptible hosts can move to the A-compartment or to the I-compartment, aware hosts can move to the S-compartment due to awareness decay or to the I-compartment due to infection (albeit at a lower rate than susceptible hosts), and infectious hosts will move to the S-compartment or A-compartment upon recovery.

The proportions of hosts in the S-, A-, and I- compartments will be denoted by $s, a, i$ respectively.
The rates of change of these fractions are governed by the following ODE model:

\begin{equation}\label{rSAISa}
\begin{split}
\frac{da}{dt} &=  \alpha_i(i) \, s \, i + \alpha_a(i) \,s \, a  + p(i) \, \delta \, i -  \beta_a a \, i - \delta_a(i) \, a, \\
\frac{di}{dt} &=  (\beta \, s  + \beta_a\, a - \delta) \, i, \qquad s+a+i = 1.
\end{split}
\end{equation}

Same assumptions about the functions $\alpha_i(i)$, $\delta_a(i)$, $\alpha_a(i)$ and the constants $\beta, \beta_a, \delta$ are made as in \cite{WJ}:  $\alpha_i(i)   \geq 0$ and $\delta_a(i) > 0$ are  differentiable functions, $\alpha_a(i) > 0$ is Lipschitz-continuous, and
$\beta, \beta_a, \delta$ are constants such that $0 \leq \beta_a < \beta$ and $0 < \delta$.
Moreover, $p(i)  \geq 0$ is a differentiable function, and $p(i) > 0$ for all $0 < i \leq 1$.

Except for $p(i)$, all these rate parameter functions and constants retain the same meanings as in \cite{WJ}: 

The term  $\alpha_i(i)i$ represents the rate at which a susceptible host becomes aware due to \emph{direct information} about the disease prevalence. 

Similarly, the term $\alpha_a(i)a$ represents the rate at which susceptible hosts become aware due to a contact with an aware host during which the latter transmits information about the disease. 

The term $\delta_a(i)$ represents the decay of awareness. It could be  a constant or any other positive Lipschitz-continuous function.

See \cite{WJ} for a discussion of how $\alpha_i(i)$, $\alpha_a(i)$ and $\delta_a(i)$ might depend on the prevalence $i$.

The inequality $\beta_a < \beta$ embodies the assumption that awareness will lead to adoption of a behavioral response that decreases the rate at which hosts contract the infection.

Finally, the term $p(i)\delta i$ embodies the \emph{direct experience} assumption, that is, a certain proportion $0 < p(i) \leq 1$ of the hosts will become aware upon recovery.  It seems plausible to assume that $p(i)$ is an increasing function in $i$. 

{\lem The region $\Omega = \{ (a,i) \in \mathbb{R}^2 \, | \, 0 \le \,  a + i \le 1, a \in [0,1] \}$ is forward-invariant. \label{invariance}}

\noindent \textit{Proof.}  By direct inspection of the system we see that $\left.(da/dt)\right|_{a=0} = \alpha_i(i) (1-i) i + p(i)\delta i \ge 0$ for $0 \le i \le 1$, $\left. (di/dt)\right|_{i=0} = 0$, and $\left. (d(a+i)/dt)\right|_{a+i=1} = p(i)\delta i - \delta_a(i)(1-i) - \delta i \leq \delta i -\delta_a(i)(1-i) - \delta i \leq 0$.
$\Box$

\subsection{Nullclines and equilibria}

To get the $i$- and $a$-nullclines, we solve the equations $di/dt = 0$ and $da/dt = 0$ respectively.

There are two  parts of the $i$-nullcline: the horizontal axis $i=0$ as the first part and the straight line
$$
i(a) = 1 -\frac{\delta}{\beta} - \left( 1 -\frac{\beta_a}{\beta} \right) a,
$$
with a slope between $-1$ and 0 under the assumption $\beta_a < \beta$ as the other part.  The expressions of the $i$-nullclines are exactly the same as in \cite{WJ}, because the second equation of \eqref{rSAISa} remain exactly the same as in the SAIS models in \cite{WJ} and not affected by $p(i)$.  They intersect at the point
$$a = \frac{\beta - \delta}{\beta - \beta_a}.$$

On the other hand, the $a$-nullcline is defined by the following equation in the variables $i$ and $a$:
$$
a^2 - \left( 1 -  i - \frac{(\alpha_i(i) + \beta_a) i + \delta_a(i)}{\alpha_a(i)} \right) a - \frac{\alpha_i(i)}{\alpha_a(i)} \,i \, \left(1-i\right) - \frac{p(i)\delta}{\alpha_a(i)}i = 0,
$$
where only the last term differs from the expression in \cite{WJ}.  Here we use the assumption that  $\alpha_a(i) > 0$.  The point $(0,0)$ is always a solution to this equation, so it is always a part of the $a$-nullcline, while the other part of the $a$-nullcline  is given by the graph of the following function~$a(i)$ on~$[0,1]$:

\begin{eqnarray*}
a(i) & = & \frac{1}{2} \, \left( 1 - i - \frac{(\alpha_i(i) + \beta_a) i + \delta_a(i)}{\alpha_a(i)} \right. \\
& & \left.
+ \   \sqrt{\left(  1 - i - \frac{(\alpha_i(i) + \beta_a) i + \delta_a(i)}{\alpha_a(i)} \right)^2 + 4 \frac{\alpha_i(i)}{\alpha_a(i)} \, i \, (1-i) + \frac{4p(i) \delta}{\alpha_a(i)}i } \right).
\end{eqnarray*}

Note that  $a(1)= \frac{1}{2} \left( -\frac{\alpha_i(1) + \beta_a + \delta_a(1)}{\alpha_a(1)} + \sqrt{\left( \frac{\alpha_i(1) + \beta_a + \delta_a(1)}{\alpha_a(1)} \right)^2 + \frac{4 p(1) \delta}{\alpha_a(1)}} \right) >0$ and $a(0) = 1 - \delta_a(0)/\alpha_a(0)$ if $\delta_a(0) \le \alpha_a(0)$, while $a(0) = 0$ if $\delta_a(0) \ge \alpha_a(0)$.  Moreover, $a(i)$ is continuous and takes on positive values for all~$i \in (0,1)$.

By sketching  these nullclines in the $i$-$a$ plane, one can see that the system has three possible types of equilibria in the first quadrant, namely,
$$
P_1=(0,0), \quad P_2 = \left( 1-\frac{\delta_a(0)}{\alpha_a(0)}, 0 \right), \quad P_3 = (a^*, i^*).
$$

\begin{figure}[H]
\begin{center}
\begin{tabular}{c}
\includegraphics[scale=0.3]{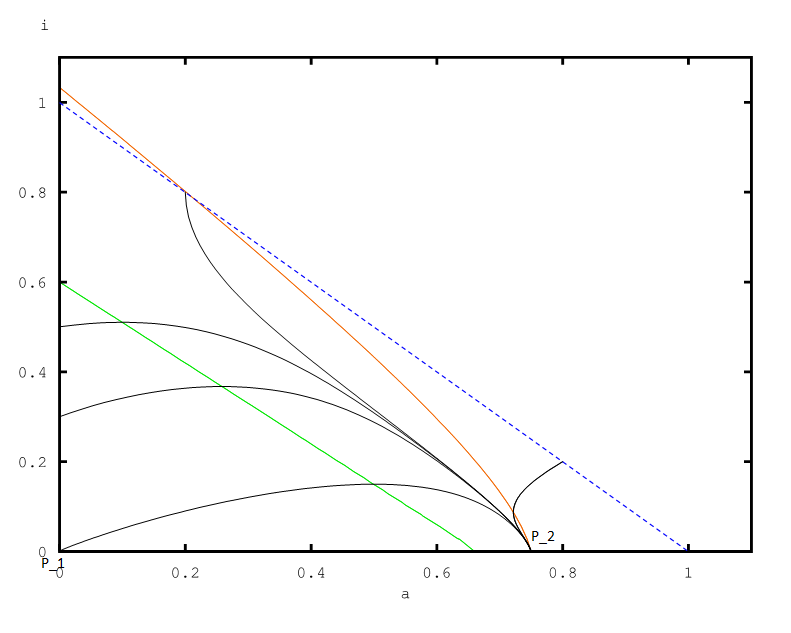}
\\
\includegraphics[scale=0.3]{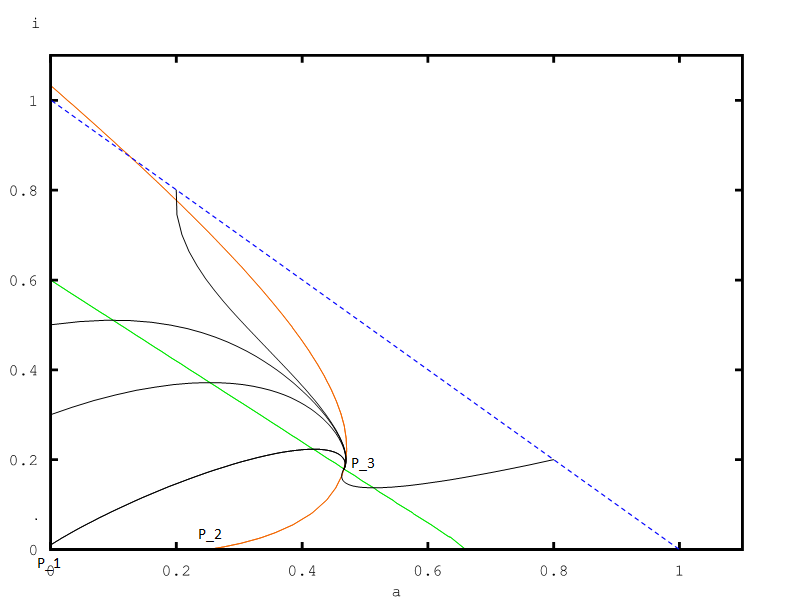}
\\
\includegraphics[scale=0.3]{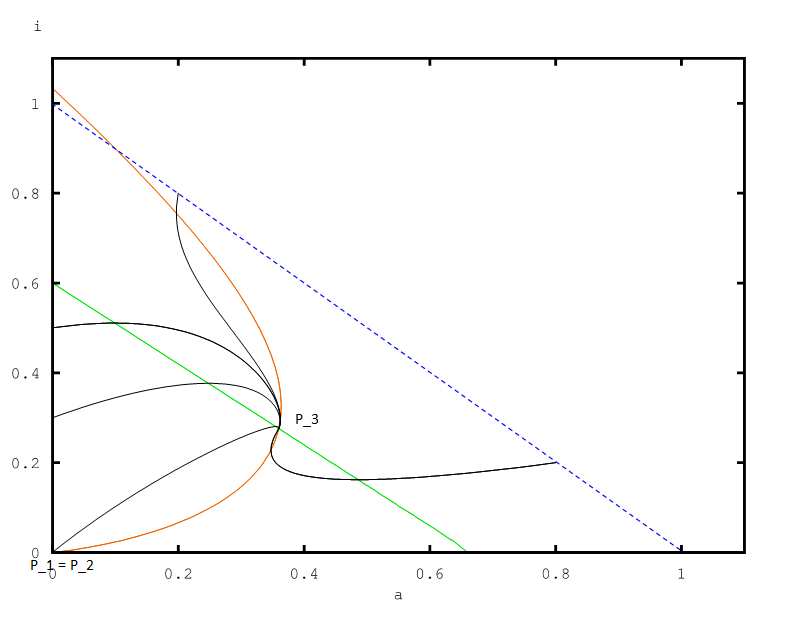}
\end{tabular}
\caption{Phase portrait of the reactive SAIAS with alert rate $\alpha_a=\alpha^0_a \, (i+1)$, $\alpha_i = \alpha_i^0 \, (i+1)$, decay rate $\delta_a=\delta_a^0/(1+i)$ and returning to awareness rate $p(i) = p_0(1+i)$ for different values of $\delta_a^0$ showing the three possible configurations of equilibria when $R_0 >1$ (top: $\delta_a^0=1$, middle: $\delta_a^0=3$, bottom: $\delta_a^0=5$). Parameters: $\delta=4$, $\beta=10$, $\beta_a=1$, $p_0 = 0.05$ and $\alpha_a^0=4$, $\alpha_i^0 = 6$.  Note that $\alpha_a(0)  > 0$  allows the existence of a second equilibrium on the $a$-axis for small values of $\delta_a^0$, which is the case in the top two panels.
\label{Portraits-rSAIS}}
\end{center}
\end{figure}

Thus, the disease-free but not awareness-free equilibrium $P_2$ exists if and only if $0~ \leq ~\frac{\delta_a(0)}{\alpha_a(0)} ~< ~1$ (see Figure \ref{Portraits-rSAIS}, in which the top two panels show parameter settings where $P_2$ exists, and the last panel shows a parameter setting without $P_2$ or where $P_1$ and $P_2$ coincide.)

Here $(a^*, i^*)$ denotes an equilibrium inside $\Omega$. Under the general assumptions made here it may not be unique, but at least one such equilibrium exists when

\begin{equation}\label{totheleft}
\frac{\beta-\delta}{\beta-\beta_a} > 1 - \frac{\delta_a(0)}{\alpha_a(0)}.
\end{equation}
Note that this condition guarantees the existence of the endemic equilibrium because the function~$a(i)$ is continuous  and satisfies
$a(1) > 0$, whereas part of the $i$-nullcline is a straight line such that $i(0) < 1$ with a slope larger than $-1$  and $a$-intercept $\frac{\beta-\delta}{\beta-\beta_a}$ (see Figure \ref{Portraits-rSAIS}, in which the top panel shows a parameter setting without $P_3$, whereas the other two panels show settings where $P_3$ exists.)

The basic reproduction numbers of the disease and awareness can be defined in the same way as in \cite{WJ}:

$$
R_0 := \beta/\delta \qquad \mbox{and} \qquad R_0^{\,a} := \alpha_a(0)/\delta_a(0),
$$
so that we can interpret the conditions for the existence of these equilibria in an intuitive way. 
The disease-and-awareness-free equilibrium $P_1$ always exists. The disease-free but not awareness-free equilibrium $P_2$ exists if and only if $R_0^{\,a} > 1$, meaning that with one aware host in a large and otherwise susceptible population,  awareness will on average increase in early stages.  The existence of an endemic equilibrium $P_3$ is guaranteed if $R_0 > \max\{1, R_0^a\}$, meaning that in early stages, the disease spreads faster than awareness, and with one infectious host in a large and otherwise susceptible population, the proportion of infectious hosts will on average increase.  However, $R_0^{\,a} > R_0 > 1$, the existence of an interior equilibrium $P_3$ is still possible as long as $\beta_a$ is close enough to $\beta$, that is, when the influence of awareness in terms of reducing the transmission rate is small enough.

The Jacobian  matrix of system \eqref{rSAISa} is
{\small
\begin{equation*}
J = \left( \begin{array}{cc}
\alpha_a(i) (s\! -\!  a) \! - \! \alpha_i(i)  i \! - \!\beta_a i \! -\! \delta_a(i)\ \ & \ \
URC
\\
- (\beta - \beta_a)  i  & \beta  s + \beta_a  a - \beta  i - \delta
\end{array}
\right),
\end{equation*}
}
with $s=1-a-i$, and the upper-right corner $URC= \alpha_i'(i) i + \alpha_i(i) - a \alpha_i'(i) i - a \alpha_i(i) - \alpha_i'(i) i^2 - 2 \alpha_i(i) i + \alpha_a'(i) a - \alpha_a'(i) a^2 - \alpha_a'(i) a i - \alpha_a(i) a  - \beta_a a - \delta_a'(i) a + p'(i) \delta i + p(i) \delta$, which is exactly the same as in the SAIS models in \cite{WJ} other than the last two terms involving~$p(i)$. 

Since $i^*=0$ at both $P_1$ and $P_2$, we get the same eigenvalues of $J$ at these equilibria as for the SAIS models in \cite{WJ}.  Specifically, the eigenvalues of $J$  at $P_1$ are
\begin{equation}\label{ewp1}
\lambda_1(P_1) = \alpha_a(0) - \delta_a(0)\qquad  \mbox{and} \qquad \lambda_2(P_1) = \beta - \delta,
\end{equation}
and the eigenvalues of $J$ at $P_2$ are
\begin{equation}\label{ewp2}
\lambda_1(P_2) =  \delta_a(0) - \alpha_a(0) \qquad \mbox{and}
\qquad
\lambda_2(P_2) = \beta - \delta - (\beta - \beta_a) \left(1 - \frac{\delta_a(0)}{\alpha_a(0)} \right).
\end{equation}

Therefore, we get the same observations as in \cite{WJ} as well:  Conditions for the existence of $P_2$ and $P_3$ imply instability of $P_1$, and condition \eqref{totheleft} implies that $\lambda_2(P_2) > 0$.  Thus, with $\lambda_1(P_2) = -\lambda_1(P_1)$, if $R_0^a > 1$, then $P_1$ is unstable and $P_2$ attracts every trajectory on the $a$-axis.

\subsection{Dynamics}

The following lemma shows that in contrast to the SAUIUAS models that we will study in Section~3, sustained oscillations are ruled out in SAIAS models.

{\lem The system \eqref{rSAISa} has no closed orbits inside $\Omega$. \label{NCO}}

\noindent \textit{Proof.}    Let $f_1(a,i)$ and $f_2(a,i)$ denote the functions on the  right-hand side of the system. The vector field $\displaystyle (F_1(a,i), F_2(a,i)) = \left( \frac{1}{a\,i}  f_1(a,i), \frac{1}{a\,i}  f_2(a,i) \right)$ is $C^1$ in the interior of $\Omega$, and
$$
\frac{\partial}{\partial a} F_1(a,i) + \frac{\partial}{\partial i} F_2(a,i) = -\frac{\alpha_i(i)}{a^2} (1-i) - \frac{\alpha_a(i)}{i} - \frac{p(i) \delta}{a^2} - \frac{\beta}{a}   <  0
$$
for all $(a,i)$ in the interior of $\Omega$.  Thus,  by Dulac's criterion of nonexistence of periodic orbits \cite{Perko}, the system \eqref{rSAISa} has no closed orbits inside $\Omega$.
$\Box$

\medskip

Since the expressions \eqref{ewp1} and \eqref{ewp2} for the  eigenvalues of the Jacobian at  $P_1$ and $P_2$ do not depend on $p(i)$, based on Lemma~2.2 and the Cantor-Bendixson Theorem, one can derive the exact analogues of behavior of trajectories in SAIAS models as for SAIS models in \cite{WJ}.


\section{SAUIUAS models}\label{SAUIS}

\subsection{The model}

Like the SAIS models of \cite{WJ}, SAIAS models ignore the degradation of information quality as it is transmitted from one individual to another.

Here we will investigate an analogy to the reactive SAIAS models of Section~\ref{SAIS-Non}, by including a compartment U of ``unwilling" hosts.  It is a generalization of W. Just and J. Salda\~{n}a's SAUIS models in \cite{WJ}, where we assume that a certain proportion of the hosts will become aware or unwilling upon recovery, instead of going back to susceptible directly.  We will name them the SAUIUAS models which are defined as follows:

\begin{equation}\label{eqn:SAUIS}
\begin{split}
\frac{da}{dt} & =  \alpha_i \, s \, i + \alpha_a \, s \, a  + p \, \delta \, i -  \beta_a \, a \, i - \delta_a\, a, \\
\frac{du}{dt} & =  \delta_a\, a  + \alpha_u \, s \, a  + q \, \delta \, i -  \beta_u\, u \, i - \delta_u\, u,\\
\frac{di}{dt} & =  (\beta  \, s + \beta_a \, a + \beta_u \, u - \delta)\, i, \qquad s+a+u+i = 1.
\end{split}
\end{equation}

Here $\alpha_u\, a$ is the rate at which susceptible hosts become unwilling after having a contact with an aware host, $\delta_u$ is the rate of awareness decay of the unwilling hosts, $p$ is the proportion of hosts that will be aware upon recovery and $q$ is the proportion of hosts that will become unwilling upon recovery, where $0 \leq p, q \leq p+q \leq 1$.  This model implicitly assumes that aware hosts  first turn into unwilling hosts before possibly entering the susceptible compartment.  Note that the SAUIS models in \cite{WJ} is a special case of the SAUIUAS models with $p = q = 0$.
All other terms play the same role as the corresponding terms in the reactive SAIAS models.

We assume that all rate constants with the possible exception of~$\beta_a, \beta_u$ are positive, and    $0 \leq \beta_a , \beta_u < \beta$.  Similarly to the SAIAS models, one could allow some of the rate coefficients to depend on~$i$.  However, in analogy with \cite{WJ}, we restrict our attention here to the case of constant rate coefficients.

{\lem The region $\Omega = \{ (a,u,i) \in \mathbb{R}^3_+ \, | \, 0 \le \,  a + u + i \le 1\}$ is positively invariant.  \label{invariance2}}

\noindent \textit{Proof.}
In the system \eqref{eqn:SAUIS}, $\left.(da/dt)\right|_{a=0} \geq 0$  and the inequality is strict when $si > 0$.  Similarly,  $\left. (du/dt)\right|_{u=0} \geq 0$ and the inequality is strict when $a > 0$. On the other hand,  $\left. (di/dt)\right|_{i=0} = 0$. It follows that the $(a,i)$- and $(u,i)$-coordinate planes  repel the trajectories  and that the $(a,u)$-plane is invariant. Now it is left to show that trajectories cannot cross the boundary $a+u+i=1$. Let $\vec{v}$ be the vector field defined by the right-hand side of the system, and let $\vec{n}=(1,1,1)$. Then with $s = 0$ in the region of the boundary $a+u+i=1$, we have  $\vec{v} \cdot \vec{n} = - \delta_u \, u - \delta \, i (1 - p - q) \leq 0$.  Therefore  the vector field on the boundary $a+u+i=1$  never points towards the exterior of $\Omega$.
$\Box$

\subsection{Possible equilibria}

System~\eqref{eqn:SAUIS} can have up to three types of equilibria in $\Omega$, namely the disease-free-and-awareness-free equilibrium $P_1 = (0,0,0)$, the disease-free equilibrium $P_2 = (a^*_0, u^*_0, 0)$, where $s^*_0 = \delta_a/\alpha_a$,
\begin{equation}\label{eqn:P2-location}
a^*_0 = \frac{ \delta_u \left( 1 - \frac{\delta_a}{\alpha_a} \right) }{
\delta_a \left( 1 + \frac{\alpha_u}{\alpha_a} \right) + \delta_u} ,
\quad
u^*_0 =
\left( 1 - \frac{\delta_a}{\alpha_a} \right)  \frac{\delta_a \left( 1 + \frac{\alpha_u}{\alpha_a} \right)}
{\delta_a \left( 1 + \frac{\alpha_u}{\alpha_a} \right) + \delta_u},
\end{equation}
and the endemic equilibrium, i.e., an interior point $P_3=(a^*,u^*,i^*)$  of $\Omega$ with

\begin{equation} \label{i*}
i^* = 1 - \left(1-\frac{\beta_a}{\beta} \right) a^* - \left(1 - \frac{\beta_u}{\beta} \right) u^* - \frac{\delta}{\beta}.
\end{equation}

Here we see that the new terms $p\delta i$ and $q\delta i$ do not affect the expressions of $a_0^*$ and $u_0^*$ relative to those in \cite{WJ}.  This is because at $P_2$, we have $i_0^* = 0$, hence $p\delta i_0^* = q\delta i_0^* = 0$.

\begin{figure}[H]
\begin{center}
\begin{tabular}{cc}
\hspace{-1cm}
\includegraphics[scale=0.6]{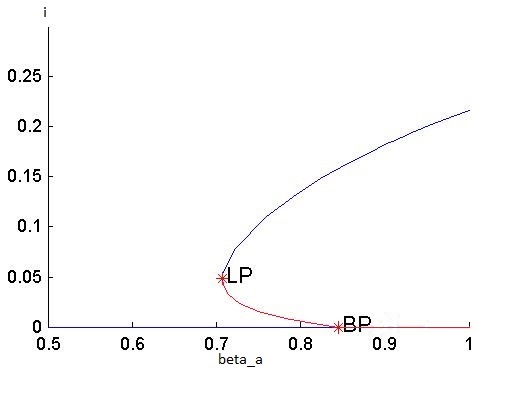}
&
\hspace{-0.75cm}
\includegraphics[scale=0.6]{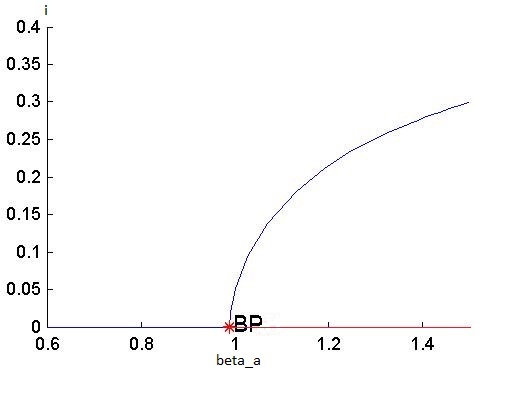}
\end{tabular}
\caption{Transcritical bifurcation diagrams of system \eqref{eqn:SAUIS} for $\beta=2$, $\beta_u=1$, $\delta=1$, $\delta_a=0.01$, $\delta_u=0.05$, $\alpha_i=0.8$, $\alpha_u=0.1$, $p = 0.1$, $q = 0.3$, $\alpha_a=0.1$ (left) and $\alpha_a=1$ (right). The stable (unstable) equilibria are depicted with a dark blue (red) line. Bifurcation values: $\beta_a^c=0.8444$ and $\beta_a=0.7071$ (for the fold bifurcation) (left panel); $\beta_a^c=0.9877$ (right panel).  
\label{Fig:Transcritical}
}
\end{center}
\end{figure}

Moreover, as the left panel of Figure~\ref{Fig:Transcritical} shows, there may be more than one interior equilibrium.

\subsection{Existence and linear stability of equilibria}

\smallskip

The disease-free-and-awareness-free equilibrium $P_1 = (0,0,0)$ always exists.
Evaluating the Jacobian matrix of system~\eqref{eqn:SAUIS} at $P_1$ we have

\begin{equation}\label{JP1}
J(P_1)  =  \left( \begin{array}{ccc}
\alpha_a - \delta_a \ \ & \ \ 0 \ \ & \ \  \alpha_i + p \delta
\\
\delta_a + \alpha_u \ \ & \ \ - \delta_u \ \ & \ \ q \delta
\\
0  \ \ & \ \ 0 \ \ & \ \ \beta - \delta
\end{array}
\right) = J_0(P_1) +  \left( \begin{array}{ccc}
0 \ \ & \ \ 0 \ \ & \ \  p \delta
\\
0 \ \ & \ \ 0 \ \ & \ \ q \delta
\\
0  \ \ & \ \ 0 \ \ & \ \ 0
\end{array}
\right) ,
\end{equation}
where $J_0(P_1)$ is the Jacobian matrix of the SAUIS models in \cite{WJ} at $P_1$.
Thus, it's eigenvalues are the same as those of $J_0(P_1)$:

$$\lambda_1(P_1)=\alpha_a - \delta_a, \quad \lambda_2(P_1)= - \delta_u, \quad \lambda_3( P_1)= \beta - \delta.$$

So, as expected, $P_1$ is unstable when $\beta > \delta$, that is, when $R_0 > 1$.  Moreover, as in the reactive SAIAS models, when  $R^a_0 := \alpha_a/\delta_a > 1$, then $P_1$ is unstable independently of the sign of $\beta - \delta$ and $P_2$ becomes biologically meaningful.

\smallskip

By~\eqref{eqn:P2-location}, the equilibrium $P_2$ exists in~$\Omega \backslash \{P_1\}$ if, and only if, $R^a_0 = \alpha_a/\delta_a > 1$.
The Jacobian matrix of~\eqref{SAUIS} at $P_2$ is
\begin{equation}\label{JP2}
J(P_2)  =  
J_0(P_2) +  \left( \begin{array}{ccc}
0 \ \,  & \ \, 0 \ \, & \ \, p\delta
\\
0 \ \, & \ \, 0 \ \, & \ \, q\delta
\\
0  \ \, & \ \, 0  \ \, & \ \, 0
\end{array}
\right),
\end{equation}

where $J_0(P_2)$ is the Jacobian matrix of the SAUIS models in \cite{WJ} at $P_2$.

The eigenvalues of $J(P_2)$ are exactly the same as those of $J_0(P_2)$.  So, $J(P_2)$ has two eigenvalues that are either negative or have negative real parts as well.

The third eigenvalue $\lambda_3(P_2)=\beta - (\beta - \beta_a) \, a^*_0 - (\beta - \beta_u) \, u^*_0 - \delta$ is negative if $\beta < \delta$ and $P_2 \in \Omega$. When $\beta > \delta$ and

\begin{equation}\label{eqn:P2lambda3+}
\frac{\beta - \beta_a}{\beta - \delta} \, a^*_0 + \frac{\beta - \beta_u}{\beta - \delta} \, u^*_0 < 1,
\end{equation}
then $\lambda_3(P_2)$ will be positive. By~\eqref{eqn:P2-location}, the values of~$a_0^*, u_0^*$ do not depend on the disease transmission parameters, and it can be seen from~\eqref{eqn:P2lambda3+} that there are large regions of the parameter space where~$\lambda_3(P_2)$ is positive and large regions where it is negative while $\beta > \delta$.

However, in terms of the existence and stability of $P_3$, more complicated situations can occur.  First, for certain regions of the parameter space, the endemic equilibrium does not exist, while for the regions of the parameter  space where an endemic equilibrium does exist, we still have subcases where such equilibrium is unique and subcases where there are more than one of them.  Examples are shown in Figure 2.  Second, in the cases for which the existence of $P_3$ is guaranteed, Hopf bifurcations may occur, and the stability of the endemic equilibrium (or equilibria when there are more than one of them) can be affected.  We will discuss Hopf bifurcation in Subsection 3.5.

\subsection{Transcritical bifurcations}

\subsubsection{Transcritical bifurcation at $R^a_0=1$}

As was stated in the previous section, the eigenvalues of the Jacobian matrix of the system~\eqref{eqn:SAUIS} at $P_1$ and $P_2$ are exactly the same as their counterparts of the SAUIS models in \cite{WJ}.  Thus, the analysis of the transcritical bifurcations at $R^a_0 = 1$ also stay the same as that in \cite{WJ}.  That is, if we assume $\beta < \delta$, then when the bifurcation parameter $R^a_0$ increases past~1, the disease-and-aware-free equilibrium $P_1$ loses its stability and $P_2$ becomes biologically meaningful and locally asymptotically stable, whereas it is unstable right before it crosses into the biologically feasible region.

\subsubsection{Transcritical bifurcations at $R_0=1$ and at $\lambda_3(P_2)=0$}
In terms of the bifurcations of an endemic equilibrium $P_3$, we adopt the same notations and strategy as those used in \cite{WJ}, where the analysis is based on the standard results for the existence of a transcritical bifurcation (see the criterion that is given right after Sotomayor's Theorem in \cite{Perko}).  Once again, calculations and results are similar to those in \cite{WJ}, i.e., we get transcritical bifurcations at $R_0 = 1$ and at $\lambda_3(P_2) = 0$.

Specifically, let~$\bf f$ denote the vector defined by the right-hand side of system \eqref{eqn:SAUIS} and let ${\bf f}_{\mu}$ be the vector of partial derivatives of its components $f_i$ with respect to a bifurcation parameter $\mu$.  Let $D{\bf f}_{\mu}$ be the Jacobian matrix of ${\bf f}_{\mu}$ and let
$D^2{\bf f}({\bf y}, {\bf y})$ be the column vector with components
$\left(D^2{\bf f}({\bf y}, {\bf y}) \right)_k := \sum_{j,l} \frac{\partial^2 f_k}{\partial x_j \partial x_l} y_j y_l$,  where $\bf y$ is a vector in $\mathbb{R}^3$, $x_1=a$, $x_2=u$, and $x_3=i$.  We use ${\bf f}^{BP}_\mu, D^2{\bf f}^{BP}$ to indicate that the above vector and function are computed at the bifurcation point.

Then on the one hand, the endemic equilibrium $P_3$ can bifurcate from $P_1$ when $\alpha_a < \delta_a$ (i.e., $R_0^a <1$). In particular, since $\lambda_3(P_1)=\beta-\delta$ is a simple eigenvalue, a bifurcation occurs for $\beta = \delta$ (i.e., at $R_0 =1$).
Taking $\beta$ as  the bifurcation parameter and evaluating the Jacobian matrix $J(P_1)$ at the bifurcation point, it follows that the row vector ${\bf u} = (0,0,1)$
 and the column vector
${\bf v} = (v_1, v_2, v_3) = (1, (\delta_a + \alpha_u + \frac{q\delta(\delta_a-\alpha_a)}{\alpha_i+p\delta})/\delta_u, (\delta_a - \alpha_a)/(\alpha_i + p\delta))^T$ are the left and right  eigenvectors for $\lambda_3=0$, respectively, where the expression of ${\bf v}$ is similar to, but not quite the same as in \cite{WJ}.  Moreover,  ${\bf f}_\beta = (0,0, (1-a-u-i)i)^T$.  Then

\vspace{0.25cm}

\noindent
1) ${\bf u} \cdot {\bf f}_{\beta}^{BP} = 0$, \\
2) ${\bf u} \cdot (D{\bf f}_{\beta}^{BP} {\bf v}) = v_3 = \frac{\delta_a - \alpha_a}{\alpha_i + p\delta} > 0$, and \\
3) ${\bf u} \cdot \left(D^2 {\bf f}^{BP}({\bf v}, {\bf v}) \right) = -2v_3 \left( (\beta-\beta_a) v_1 + (\beta-\beta_u) v_2 + \beta v_3 \right) < 0$.

\vspace{0.25cm}
\noindent

Thus, when $R^a_0 < 1$  system \eqref{eqn:SAUIS} experiences a transcritical bifurcation as $\beta$  crosses the bifurcation value $\beta=\delta$ \cite{Perko}.
Moreover,  Theorem 4.1 in~\cite{Castillo}, together with the inequality in~2) and the inequality in~3), implies that the direction of the bifurcation is always the same, namely, system~\eqref{SAUIS} experiences a forward bifurcation at~$R_0 = 1$.

On the other hand, assume $\alpha_a > \delta_a$ such that a positive $P_2$ exists, and $\beta > \delta$ such that $\lambda_3(P_2)$ can be positive for some parameters values. From the discussion surrounding~\eqref{eqn:P2lambda3+}  it follows that $\lambda_3(P_2)=0$ if and only if
\begin{equation}\label{BP}
\frac{\beta - \beta_a}{\beta - \delta} \, a^*_0 + \frac{\beta - \beta_u}{\beta - \delta} \, u^*_0 = 1,
\end{equation}
with $a^*_0$ and $u^*_0$ given by \eqref{eqn:P2-location}. That is, at this parameter combination, the endemic equilibrium $P_3$ bifurcates from $P_2$.

\bigskip

\vspace{0.25cm}
\noindent
Applying the same calculation and argument as in \cite{WJ}, where the relevant objects are not affected by $p$ and $q$, we conclude that in general, for a nonempty open set of parameter settings at which~$P_1 \neq P_2 \in \Omega$ the  system~\eqref{eqn:SAUIS} experiences a transcritical bifurcation as $\beta_a$ passes through the bifurcation value
$$
\beta_a^c :=\beta - \frac{1}{a_0^*} \left( \beta-\delta - (\beta-\beta_u)u_0^* \right).
$$
But in contrast to what happens at $R_0=1$, the direction of the bifurcation is not always the same.  An example of forward and backward bifurcations occurring at $\lambda_3(P_2)=0$ for $\beta_a < \beta$ is shown in Figure~\ref{Fig:Transcritical}.

The same conclusion holds if we use $\beta$ or $\beta_u$ as a  bifurcation parameter.

\subsection{Hopf bifurcations}

It is shown in \cite{WJ} that Hopf bifurcations and sustained oscillations are possible in the SAUIUAS models when $p = q = 0$.  
To explore the influence of $p$ and $q$ on the occurrence of Hopf bifurcations, we adopt the same strategy as in \cite{WJ}.
We call a pair $(\sigma, \tau)$ of parameters of~\eqref{eqn:SAUIS} a \emph{Hopf pair} if there exists an equilibrium point  at which the Jacobian matrix has a pair of pure imaginary eigenvalues. 

An explicit criterion that specifies whether an $n \times n$ matrix~$M$, with coefficients that may depend upon parameters, has a pair of pure imaginary eigenvalues is given in \cite{GMS97}.

To locate the Hopf pairs in a $(\sigma, \tau)$ parameter space, we use~$\sigma$ as free parameter and solve the system given by  the equilibrium equations combined with the conditions in the criterion mentioned above in \cite{GMS97} 
for $a^*$, $u^*$, $i^*$, and $\tau$, while all other parameters are set to fixed values. 

The set of Hopf pairs defines the so-called Hopf- bifurcation curve $H$ in $(\sigma, \tau)$ parameter space.

$H$ can be parametrized by $i^*$ (or any of the components of $P_3$) \cite{SSK12}. 

\bigskip

Let us first take $(p,q)$ as the Hopf pair, while fixing all other parameters.   Two examples of curve $H$ in the $(p,q)$ parameter space are shown in Figure~ 3 and Figure~ 7 in the cases of $\alpha_u = 1$ and $\alpha_u = 3$, where $\delta = 1$, $\delta_a = 0.01$, $\delta_u = 0.05$, $\beta = 3$, $\beta_u = 0.5$, $\alpha_a = 0.012$, $\beta_a=0.2$ and $\alpha_i=0.05$.

For $\alpha_u = 1$, the prevalence of the disease at the bifurcation points along the curve $H$ is presented in Figure~ 4.

In Figure~ 3, $P_3$ is unstable under the curve and then is stable above the curve
.  For example, when $p = 0.05$, $q^*=0.1075$ is the only Hopf-bifurcation point.  Figure~ 5 shows the Hopf-bifurcation diagram with $p = 0.05$, where $P_3$ is unstable and we can get sustained oscillations when $q < q^*$ and $P_3$ is stable when $q > q^*$.  

Numerical simulations confirm our predictions.  Examples are shown in Figure~ 6, where sustained oscillations can be observed when $q = 0.05$ and 0.1 in the top panels, and $P_3$ becomes stable in the bottom panels when we increase the value of $q$ such that it exceeds $q^* = 0.1075$.  Specifically, $q = 0.15$ in the bottom left panel and 0.2 in the bottom right panel.

\medskip

For $\alpha_u = 3$, the prevalence of the disease at the bifurcation points along the curve $H$ is presented in Figure~ 8.

In Figure~ 7, $P_3$ is unstable under the curve and is stable above the curve
.  For example, when $p = 0.05$, $q^*=0.1923$ is the only Hopf-bifurcation point.  Also note that under the parameter setting of Figure~ 7, sustained oscillations can be observed even when $p =1$.  Figure~ 9 shows the Hopf-bifurcation diagram with $p = 0.05$, where $P_3$ is unstable and we can get sustained oscillations when $q < q^*$ and $P_3$ is stable when $q > q^*$.  

Numerical simulations confirm our predictions.  Examples are shown in Figure~ 10, where with $p = 0.05$ sustained oscillations can be observed when $q = 0.05$ and 0.15 in the top panels, and $P_3$ becomes stable in the bottom left panel when we increase the value of $q$ to 0.25 such that it exceeds $q^* = 0.1923$.  
With $p = 1$ and $q = 0$, we also see sustained oscillations in the bottom right panel.  

In addition, Figure~3 and Figure~7 indicate the possibility of getting Hopf bifurcation at $p = 1$ and $q = 0$ for some $1 < \alpha_u < 3$ while other parameters remain the same.  
This $\alpha_u$ is found to be  2.1015.

\bigskip

Observing Figure~ 3 and Figure~ 7, it seems that in order to get a Hopf bifurcation, the upper bound of $q$ cannot be large when $\alpha_u$ is not too large.  To see the relationship between $q$ and $\alpha_u$ when Hopf-bifurcations occur, we assume $p = 1-q$, meaning that a host becomes either aware or unwilling upon recovery, and take $(q,\alpha_u)$ as the Hopf pair, while fixing all other parameters.   An example of curve $H$ in the $(q,\alpha_u)$ parameter space is shown in Figure~ 11, where  $\delta = 1$, $\delta_a = 0.01$, $\delta_u = 0.05$, $\beta = 3$, $\beta_a = 0.2$, $\beta_u = 0.5$, $\alpha_a = 0.012$, $\alpha_i = 0.05$.

The prevalence of the disease at the bifurcation points along the curve $H$ is presented in Figure~ 12.  We can see clearly how $i^*$ decreases as $\alpha_u$ increases.

In Figure~ 11, $P_3$ is stable under the curve
, and is unstable above the curve.  When $q = 0.1$, $\alpha_u^*=2.6449$ is the only Hopf-bifurcation point. 

Numerical simulations confirm our predictions.  Examples are shown in Figure~ 13, where $P_3$ is stable when $\alpha_u = 0.5$ and 2 in the top panels, and sustained oscillations can be observed in the bottom panels when we increase the value of $\alpha_u$ such that it exceeds $\alpha_u^* = 2.6449$.  Specifically, $\alpha_u = 3$ in the bottom left panel and 5 in the bottom right panel.

\bigskip

Further, if we bound $\alpha_u$ by 5 from above, still set $\delta = 1$, $\delta_a = 0.01$, $\delta_u = 0.05$, $\beta = 3$, $\beta_a = 0.2$, $\beta_u = 0.5$, $\alpha_a = 0.012$, and $\alpha_i = 0.05$, the Hopf-curve in the parameter space $(\alpha_u, q)$ given by Figure~ 14 shows a ``zoomed in" relationship between $\alpha_u$ and $q$ under the assumption that $p = 1-q$.  Note that in order to do this, we switched the order of the parameters in the previous Hopf pair, i.e., we took $(\alpha_u,q)$ as the Hopf pair.  This is because in our numerical code, the parameter on the vertical axis is a dependent variable, whose value is calculated for each value of the chosen parameter on the horizontal axis playing the role of an independent variable.   Thus, in order to bound $\alpha_u$ from above by 5, we need $\alpha_u$ to be an independent variable and set it as the parameter on the horizontal axis.

The prevalence of the disease at the bifurcation points along the curve H is presented in Figure~ 15. We can see clearly how $i^*$ decreases as $q$ increases.

In Figure~ 14, $P_3$ is stable above the curve and is unstable under the curve.  When $\alpha_u = 3$, $q^*=0.1549$ is the only Hopf-bifurcation point. 

Examples of numerical simulations are shown in Figure~ 16, where sustained oscillations can be observed when $q = 0.05$ and 0.1 in the top panels, and $P_3$ becomes stable in the bottom panels when we increase the value of $q$ such that it exceeds $q^* = 0.1549$.  Specifically, $q = 0.5$ in the bottom left panel and 0.8 in the bottom right panel.

\bigskip

In fact, in order to explore the influence of awareness that results from direct experience under the assumption that $p = 1-q$, each $(\sigma, q)$ can be a natural choice as the Hopf pair.  For example, we take $(\alpha_i, q)$ as the Hopf pair, with  $p = 1-q$, $\delta = 1$, $\delta_a = 0.01$, $\delta_u = 0.05$, $\beta = 3$, $\beta_a = 0.2$, $\beta_u = 0.5$, $\alpha_a = 0.012$.  If $\alpha_u = 1$, there is no Hopf-bifurcation, and $P_3$ is stable.  If $\alpha_u = 3$, we get the Figures~ 17--19.

\bigskip

Moreover, now that we have shown the possibility of getting sustained oscillations and Hopf bifurcations with the awareness gained from direct experience taken into consideration, including the case of $p = 1$ and $q = 0$, we can also explore the dynamics of such models while assuming that all the infectious hosts will be unwilling upon recovery.  That is, we consider the SAUIUAS models under the assumption that $q = 1$, and see whether sustained oscillations can still be observed. 

\bigskip 

In fact, we can get sustained oscillations in the SAUIUAS models when $p = 0$ and $q = 1$.  If we set  $\beta = 3$, $\beta_a = 0.2$, $\beta_u = 0.4$, $\alpha_i = 0.05$, $\alpha_a = 0.012$, $\delta_a = 0.01$, $\delta_u = 0.05$, $\delta = 1.7$, and $\alpha_u = 30$ and take $(p, q)$ as the Hopf pair, we get the Hopf-bifurcation curve in Figure~ 20, where sustained oscillations can be observed for $(p, q)$ under the curve.  This under-the-curve region covers the entire biologically realistic region where $p \geq 0$, $q \geq 0$, and $p + q \leq 1$.   

Numerical simulations are shown in Figure 21.  Sustained oscillations occur under the above parameter settings and $p = 0$, $q = 1$.

\bigskip 

Finally fix $p = 0$ and $q = 1$, and take $(\delta, \alpha_u)$ as the Hopf pair.  We get Figure 22--25, and we can see how $\alpha_u$ decreases as $\delta$ increases when Hopf-bifurcations are observed.

\section{Discussion}
Here we showed that the SAIAS and SAUIUAS models exhibit the same rich dynamics as the SAIS and SAUIA models in \cite{WJ}.  This indicates that in order to get sustained oscillations, the unrealistic assumption that all infected hosts will return to the susceptible compartment directly upon recovery without any awareness gained from direct experience is not necessary.  For a detailed discussion of the significance of the observed patterns and the relation of these findings to the broader literature, see Section 4 of \cite{WJ}.

Moreover, increasing the proportions of hosts that become aware or unwilling upon recovery in different ways or in different parameter regions can have various effects.  For example, in the biologically feasible region of Figure 7, start from any point under the curve, increase the value of $q$ while $p$ is fixed.  Sustained oscillations will disappear and an endemic equilibrium will become stable when crossing the curve.  However, if we start from a point to the left of the curve, increase the value of $p$ while fixing $q$, then an initially stable endemic equilibrium becomes unstable and sustained oscillations are born while crossing the curve.  But the amplitudes of the oscillations will decrease when the value of $p$ is further increased.  Finally, if we start from a point above the curve, then increasing $p$ or $q$ will not change the stability of the endemic equilibrium, only move $i^*$ to a lower level. 

However, our SAIAS and SAUIUAS models are still overly simplified in many aspects.  For example, we implicitly based our models on the uniform mixing assumption.  So a possible next step is to develop and investigate their network-based versions.  Another issue is that in the SAUIUAS models, we set all the parameters as constants, while some of them are more likely to be non-constant functions of $i$.  It will be of interest to explore whether non-constant rate functions can lead to even richer dynamics.

\bigskip

\noindent
\textbf{Acknowledgements}\ \ I hereby thank Professor Winfried Just of Ohio University for helping me to formulate the questions that were explored here, carefully reading earlier drafts of this manuscript, and sharing his comments.  I also thank Professor Joan Salda\~{n}a of the University of Girona for sharing his MATLAB code with me, helping me with adapting it for the numerical explorations reported here, and sharing his comments with me after reading an earlier draft carefully.

\begin{figure}[H]\label{Fig:Hopf-p-q-1}
\begin{center}
\begin{tabular}{c}
\includegraphics[scale=0.2]{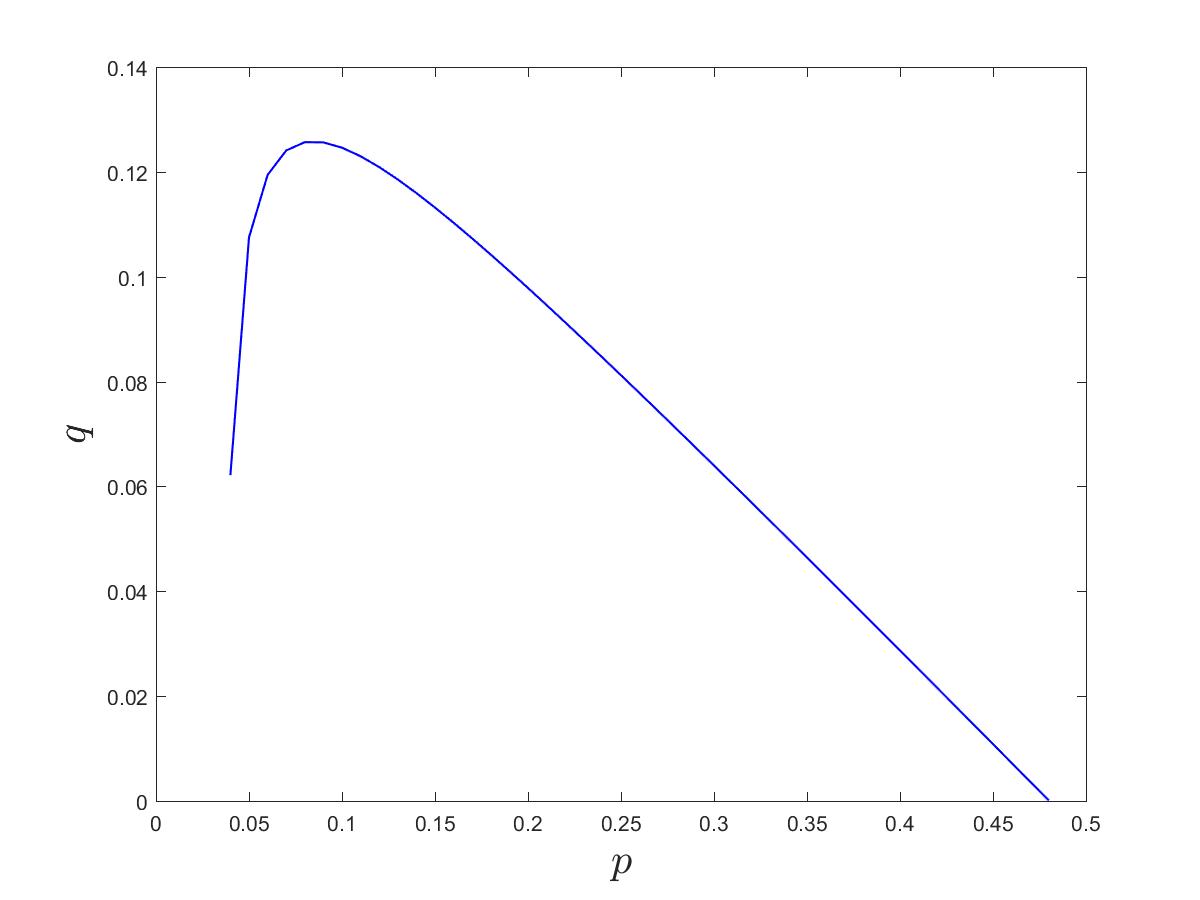}
\end{tabular}
\caption{Hopf-bifurcation curve. $p$: $0.04$(left) -- $0.48$(right).  When $p = 0.05$, $q^* = 0.1075$ is the only Hopf-bifurcation point. Here, $\alpha_u = 1$, $\delta = 1$, $\delta_a = 0.01$, $\delta_u = 0.05$, $\beta = 3$, $\beta_u = 0.5$, $\alpha_a = 0.012$, $\beta_a = 0.2$, $\alpha_i = 0.05$.}
\end{center}
\end{figure}

\begin{figure}[H]\label{Fig:Hopf-infectious-q-1}
\begin{center}
\begin{tabular}{c}
\includegraphics[scale=0.2]{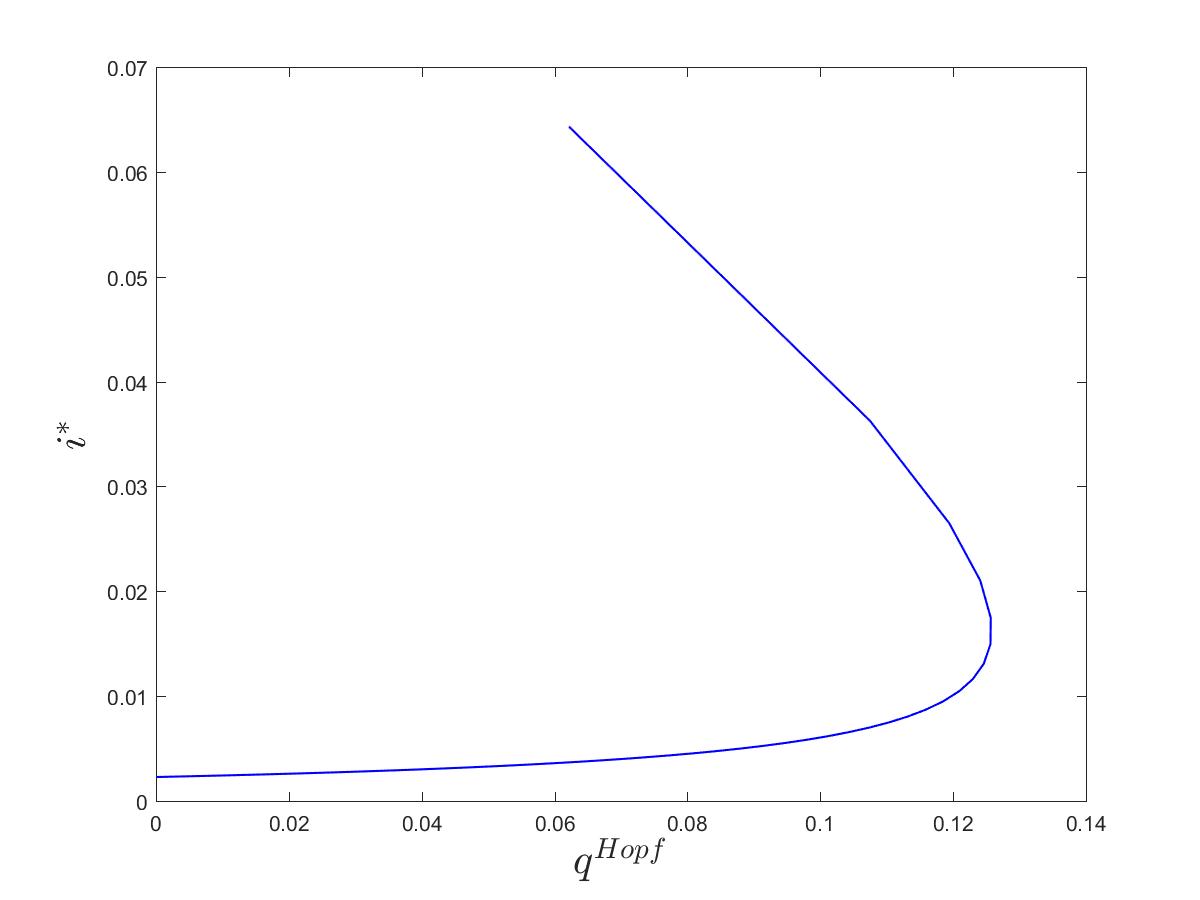}
\end{tabular}
\caption{Fraction of infectious hosts as a function of $q$ along the Hopf-bifurcation curve in Figure 3, where $\alpha_u = 1$, $\delta = 1$, $\delta_a = 0.01$, $\delta_u = 0.05$, $\beta = 3$, $\beta_u = 0.5$, $\alpha_a = 0.012$, $\beta_a = 0.2$, $\alpha_i = 0.05$. }
\end{center}
\end{figure}

\begin{figure}[H]\label{Fig:Hopf-infectious-diagram-q-1}
\begin{center}
\begin{tabular}{c}
\includegraphics[scale=0.3]{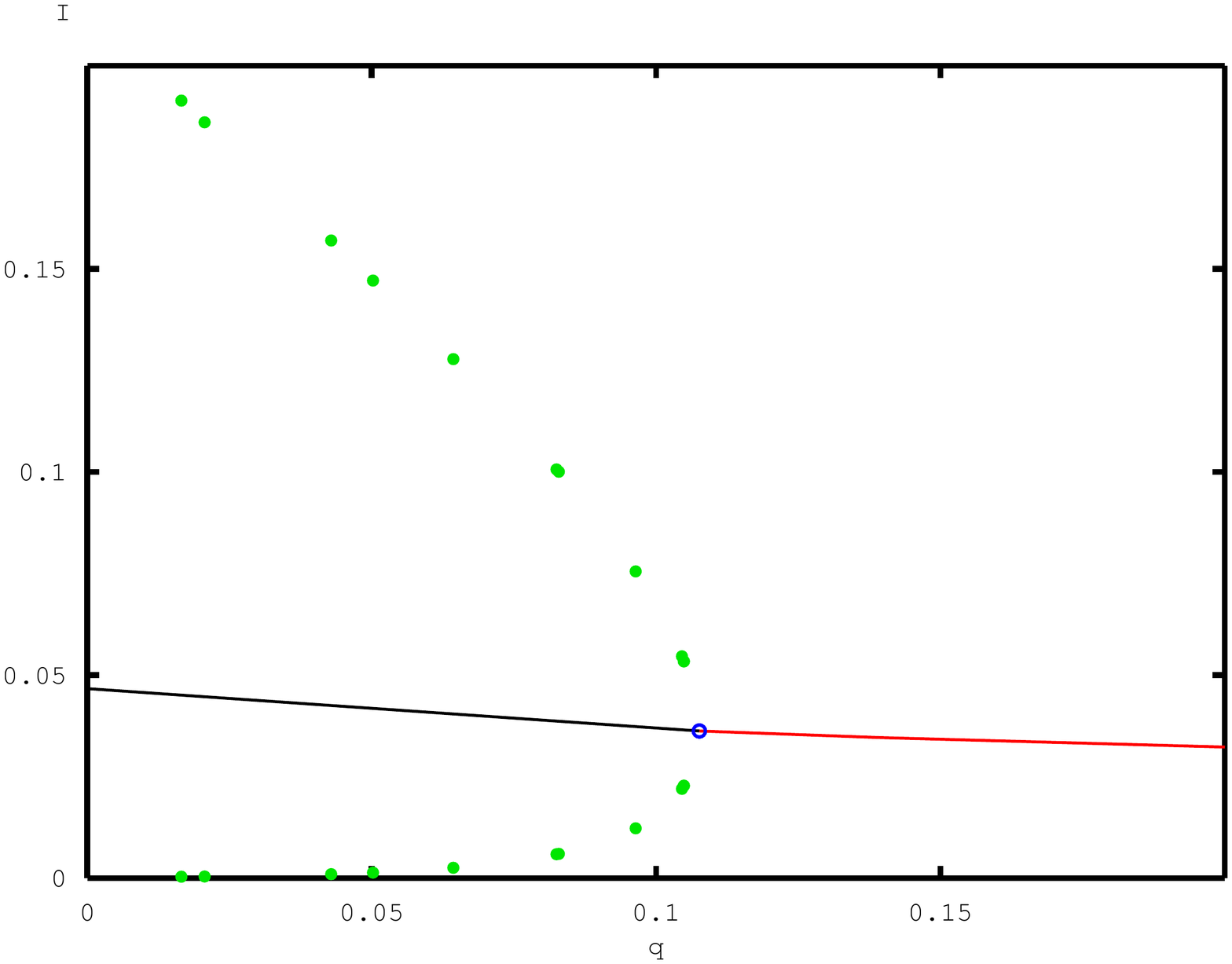}
\end{tabular}
\caption{Hopf-bifurcation diagram -- fraction of infectious hosts as a function of $q$ for $p=~0.05$, where stable endemic equilibria are depicted in red and the unstable ones are in black.  The green dots indicate the boundaries of oscillations corresponding to each value of $q$ when the corresponding endemic equilibrium is unstable and oscillations occur. Here, with $\alpha_u = 1$, $\delta = 1$, $\delta_a = 0.01$, $\delta_u = 0.05$, $\beta = 3$, $\beta_u = 0.5$, $\alpha_a = 0.012$, $\beta_a = 0.2$, $\alpha_i = 0.05$, when $p = 0.05$, $q^* = 0.1075$ is the only Hopf-bifurcation point.}
\end{center}
\end{figure}

\begin{figure}[H]
\begin{center}
\begin{tabular}{cc}
\hspace{-1cm}
\includegraphics[scale=0.2]{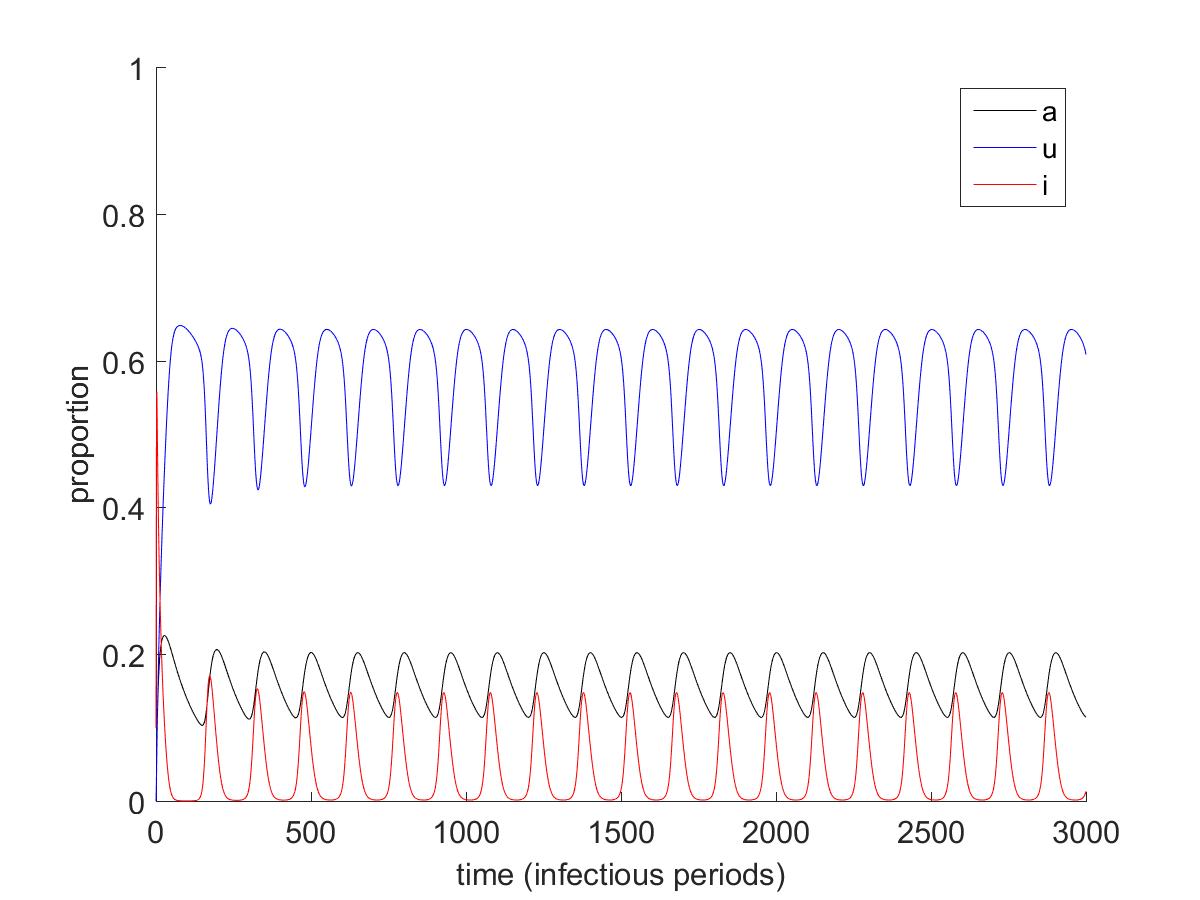}
&
\hspace{-0.75cm}
\includegraphics[scale=0.2]{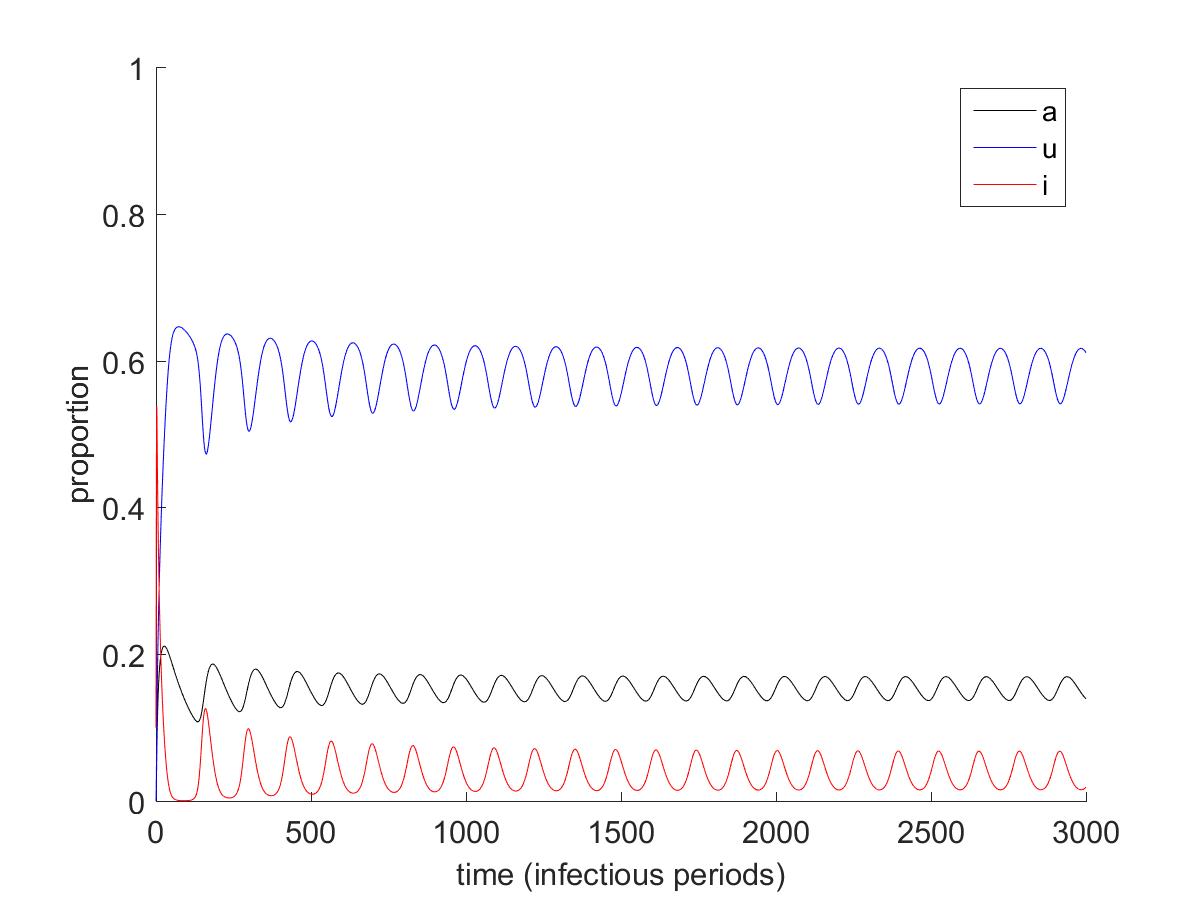}
\\
\hspace{-1cm}
\includegraphics[scale=0.2]{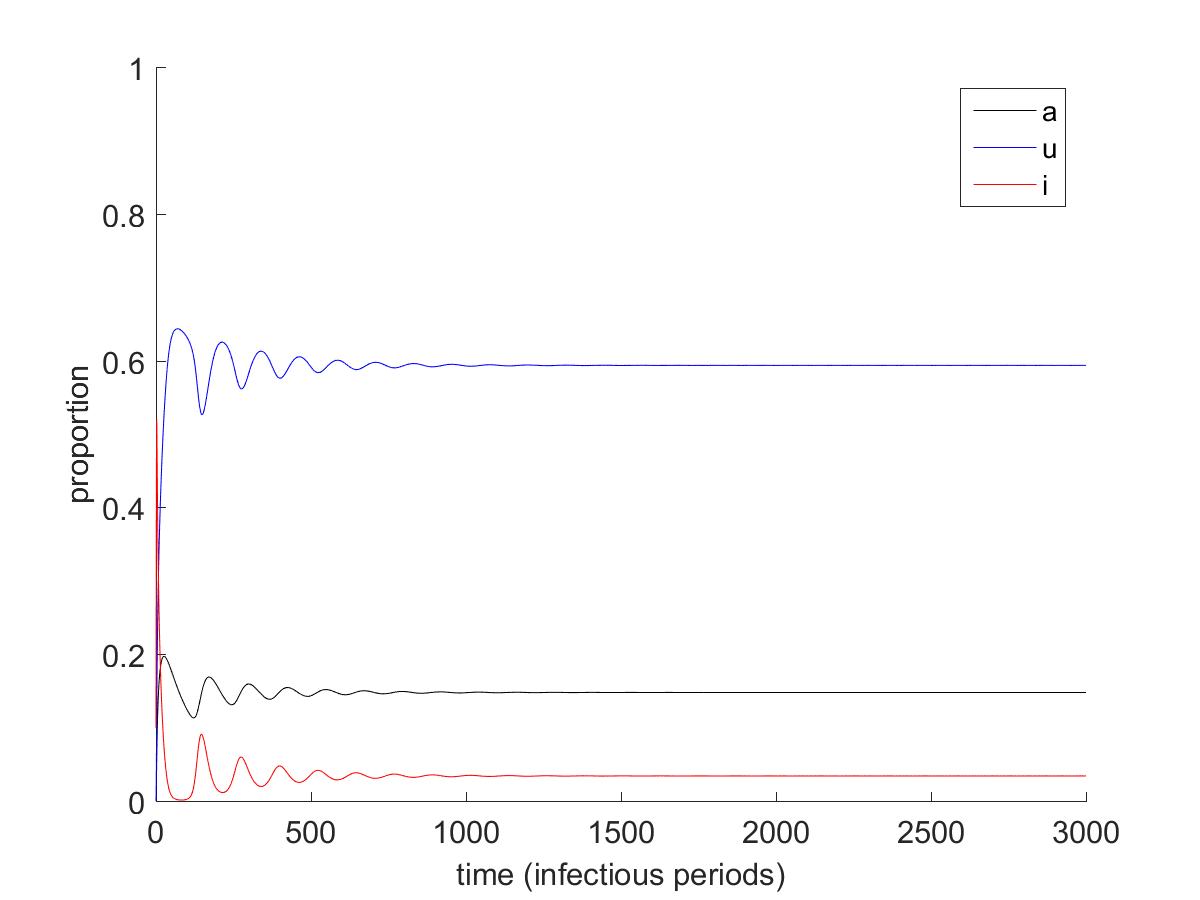}
&
\hspace{-0.75cm}
\includegraphics[scale=0.2]{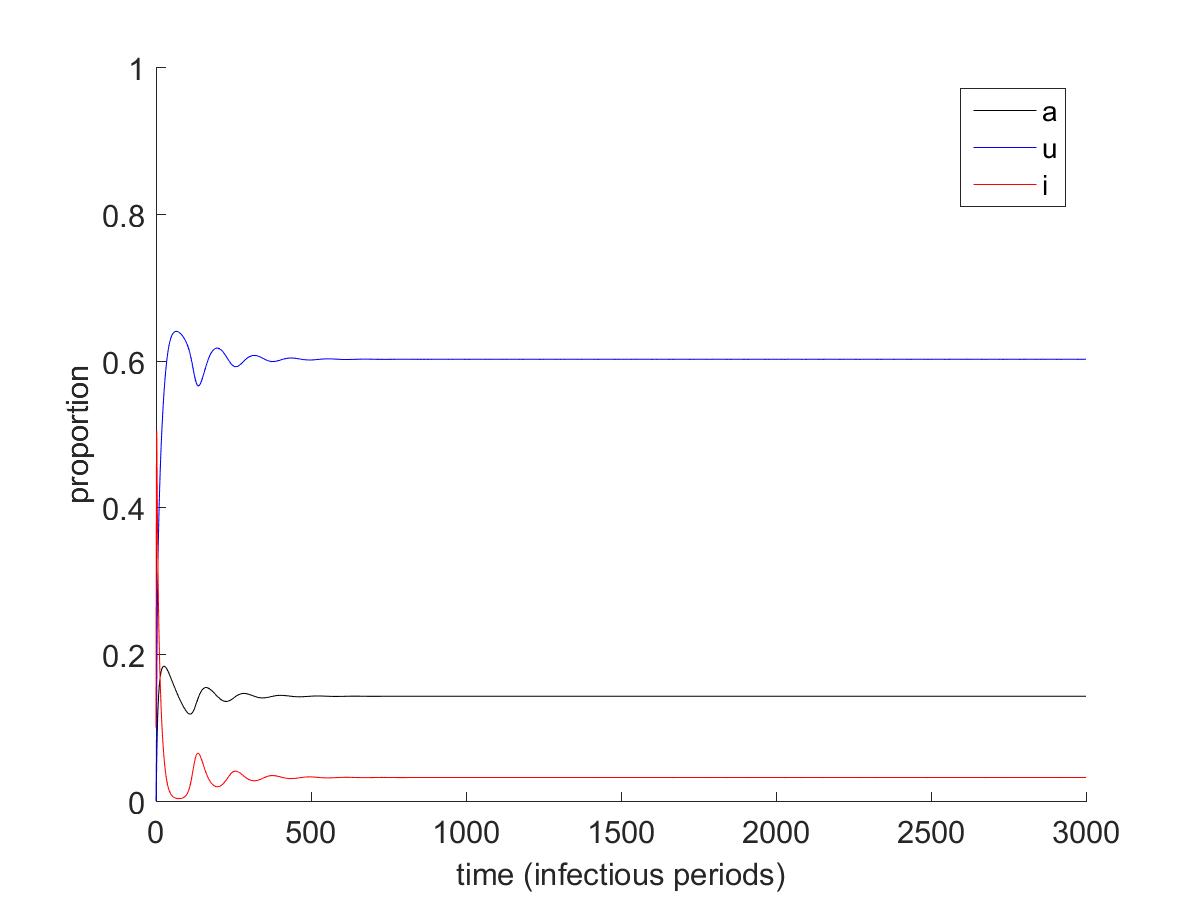}
\end{tabular}
\caption{Evolution of the fraction of infectious, aware, and unwilling hosts for different values of $q$ along the vertical section in the Hopf-bifurcation curve corresponding to $p=0.05$ (when $p = 0.05$, $q^* = 0.1075$ is the only Hopf-bifurcation point): $q=0.05$ (top left), 0.1 (top right), 0.15 (bottom left), 0.2 (bottom right).  Initial condition:  $a(0) = u(0) = 0$, $i(0) = 0.1$.  Here, $\alpha_u = 1$, $\delta = 1$, $\delta_a = 0.01$, $\delta_u = 0.05$, $\beta = 3$, $\beta_u = 0.5$, $\alpha_a = 0.012$, $\beta_a=0.2$ and $\alpha_i=0.05$.
\label{Fig:Evolution-SAUIS}}
\end{center}
\end{figure}

\begin{figure}[H]
\begin{center}
\begin{tabular}{c}
\includegraphics[scale=0.27]{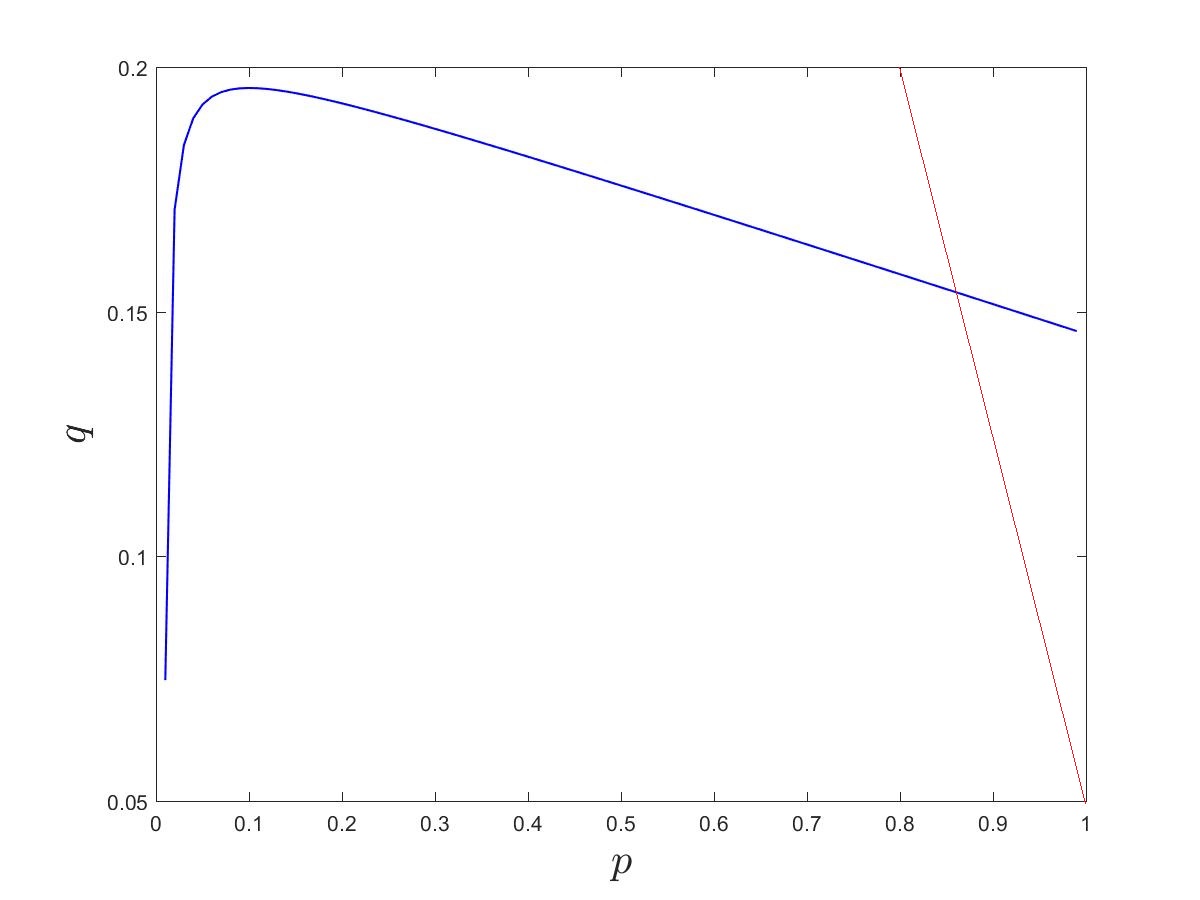}
\end{tabular}
\caption{Hopf-bifurcation curve. $p$: $0.01$(left) -- $1$(right).  The red line sets the boundary $p + q = 1$, so only the part of the Hopf-bifurcation curve that lies on the left of the red line makes biological sense.  When $p = 0.05$, $q^*=0.1923$ is the only Hopf-bifurcation point.  Here, $\alpha_u = 3$, $\delta = 1$, $\delta_a = 0.01$, $\delta_u = 0.05$, $\beta = 3$, $\beta_u = 0.5$, $\alpha_a = 0.012$, $\beta_a = 0.2$, $\alpha_i = 0.05$.
}
\end{center}
\end{figure}

\begin{figure}[H]
\begin{center}
\begin{tabular}{c}
\includegraphics[scale=0.2]{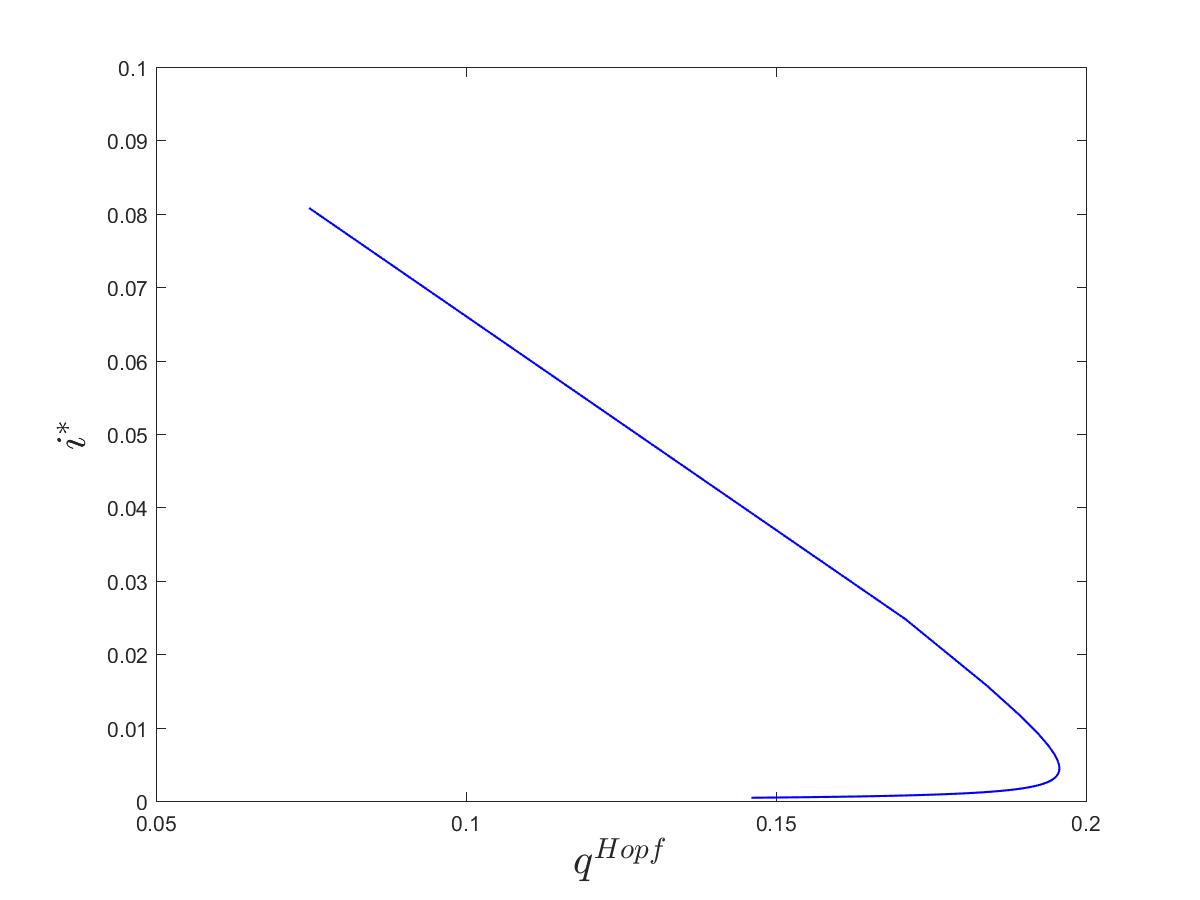}
\end{tabular}
\caption{Fraction of infectious hosts as a function of $q$ along the Hopf-bifurcation curve in Figure 7, where $\alpha_u = 3$, $\delta = 1$, $\delta_a = 0.01$, $\delta_u = 0.05$, $\beta = 3$, $\beta_u = 0.5$, $\alpha_a = 0.012$, $\beta_a = 0.2$, $\alpha_i = 0.05$.
}
\end{center}
\end{figure}

\begin{figure}[H]
\begin{center}
\begin{tabular}{c}
\includegraphics[scale=0.3]{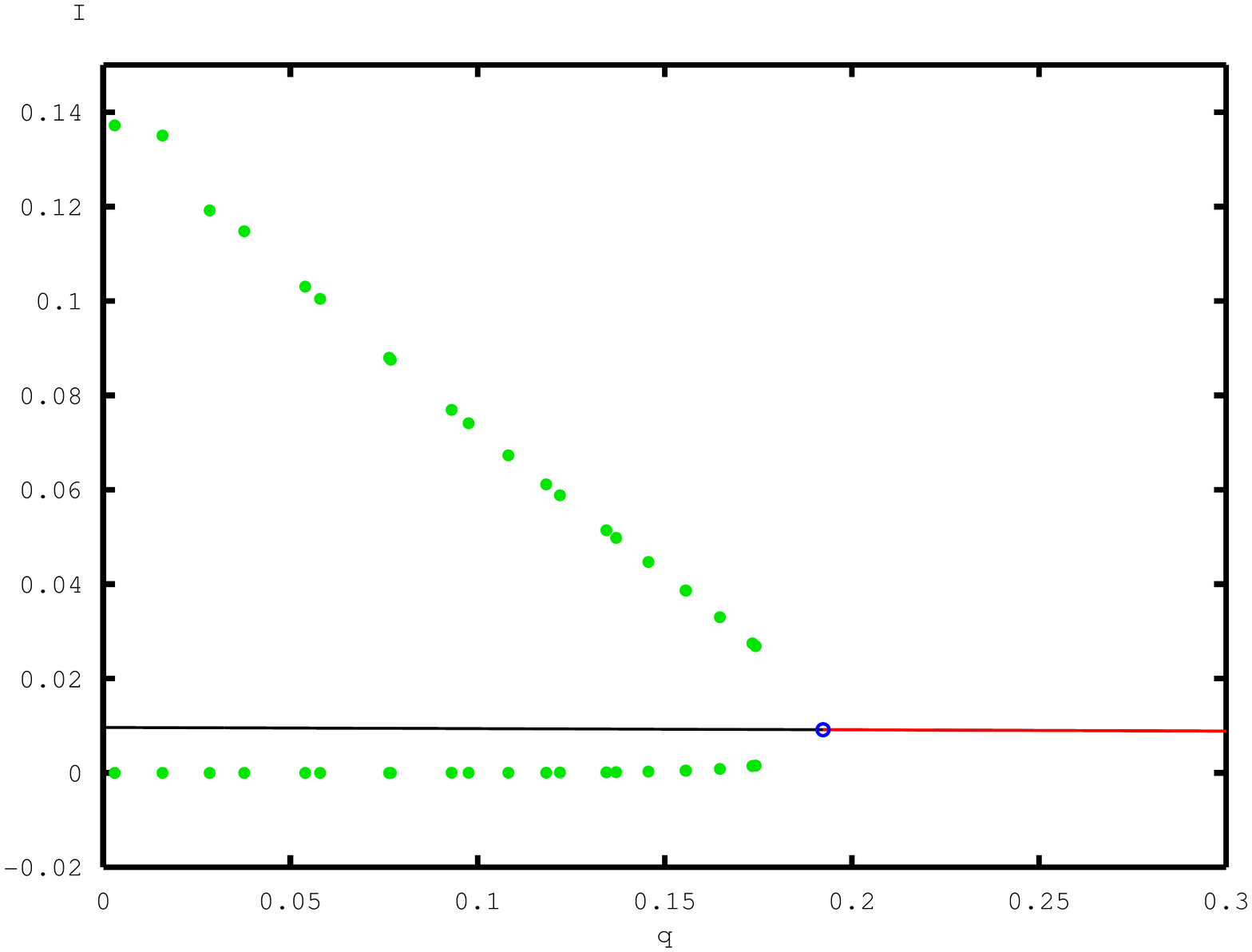}
\end{tabular}
\caption{Hopf-bifurcation diagram -- fraction of infectious hosts as a function of $q$ for $p=~0.05$, where stable endemic equilibria are depicted in red and the unstable ones are in black.  The green dots indicate the boundaries of oscillations corresponding to each value of $q$ when the corresponding endemic equilibrium is unstable and oscillations occur.  Here, with $\alpha_u = 3$, $\delta = 1$, $\delta_a = 0.01$, $\delta_u = 0.05$, $\beta = 3$, $\beta_u = 0.5$, $\alpha_a = 0.012$, $\beta_a = 0.2$, $\alpha_i = 0.05$, when $p = 0.05$, $q^* = 0.1923$ is the only Hopf-bifurcation point.
}
\end{center}
\end{figure}

\begin{figure}[H]
\begin{center}
\begin{tabular}{cc}
\hspace{-1cm}
\includegraphics[scale=0.2]{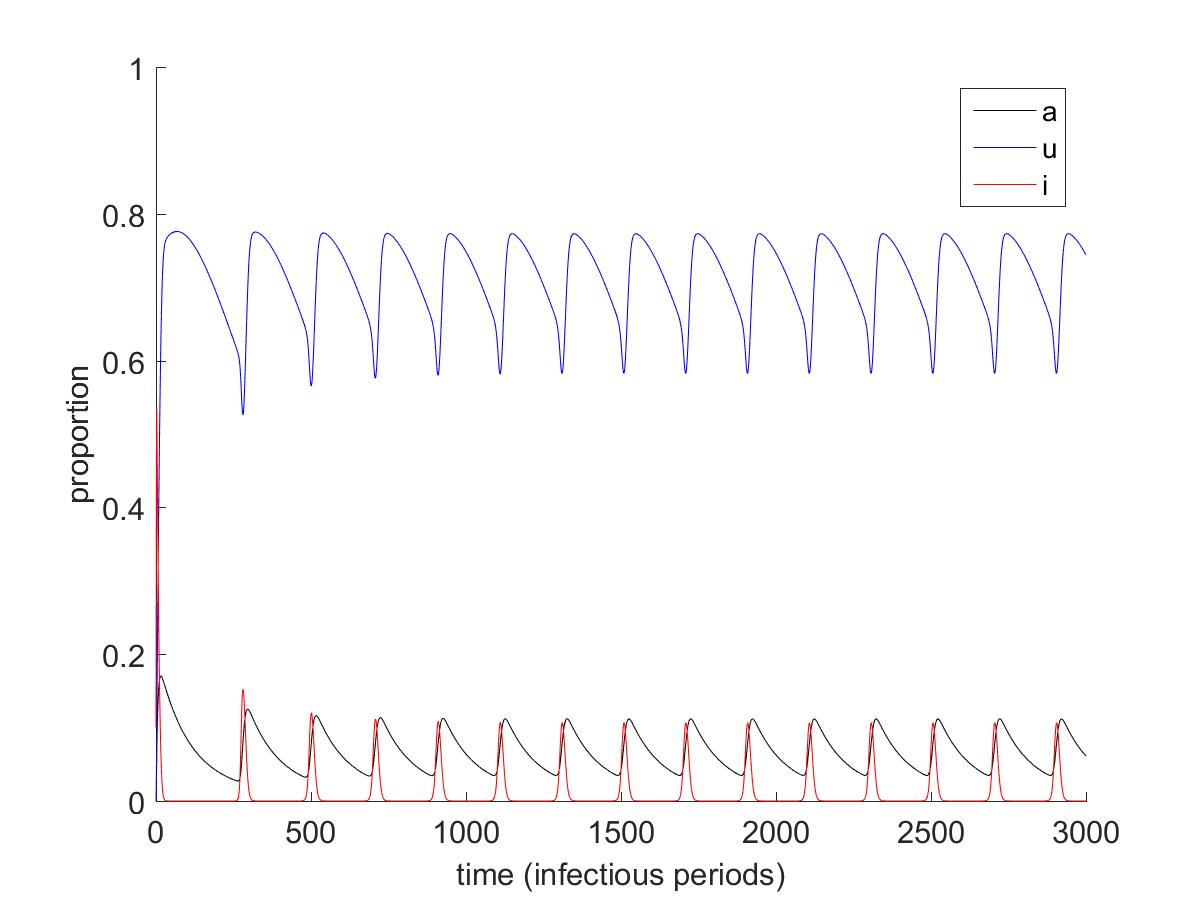}
&
\hspace{-0.75cm}
\includegraphics[scale=0.2]{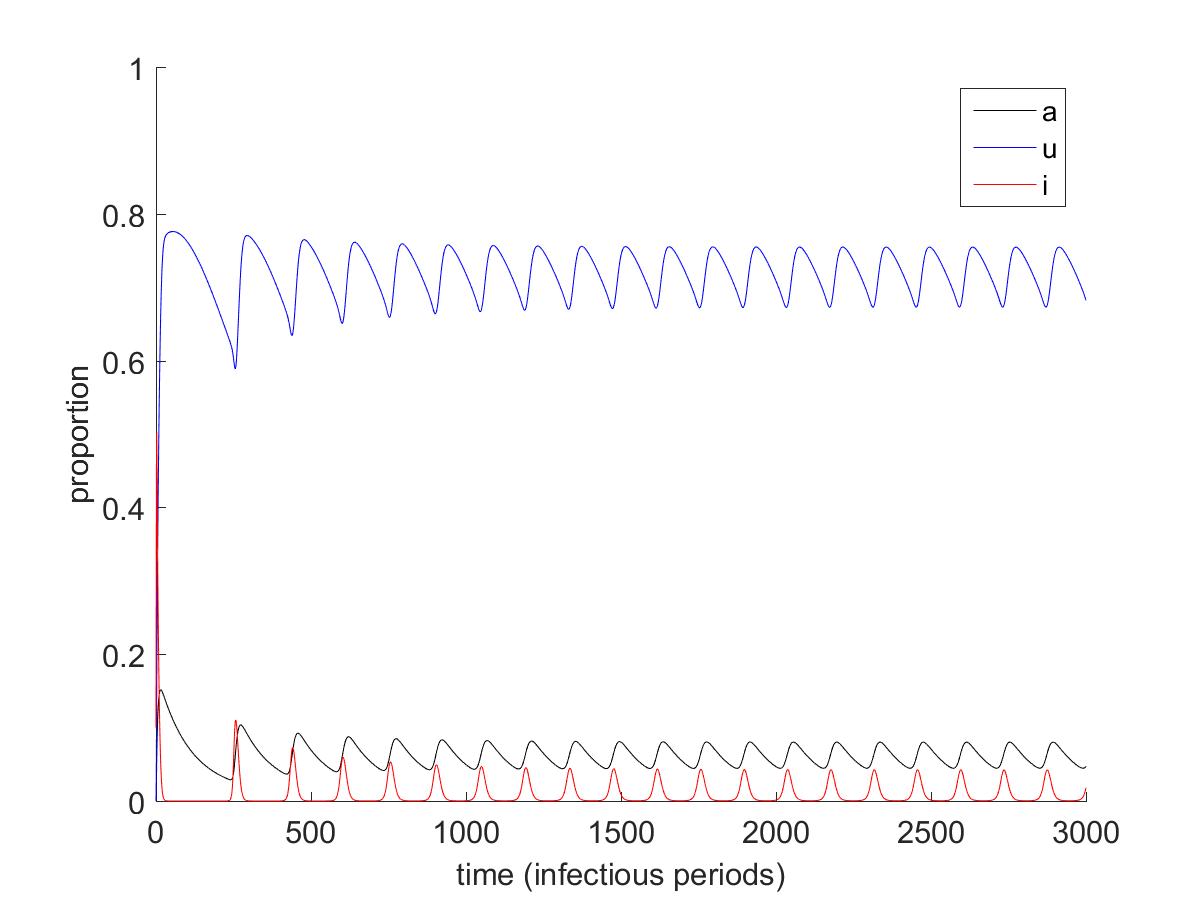}
\\
\hspace{-1cm}
\includegraphics[scale=0.2]{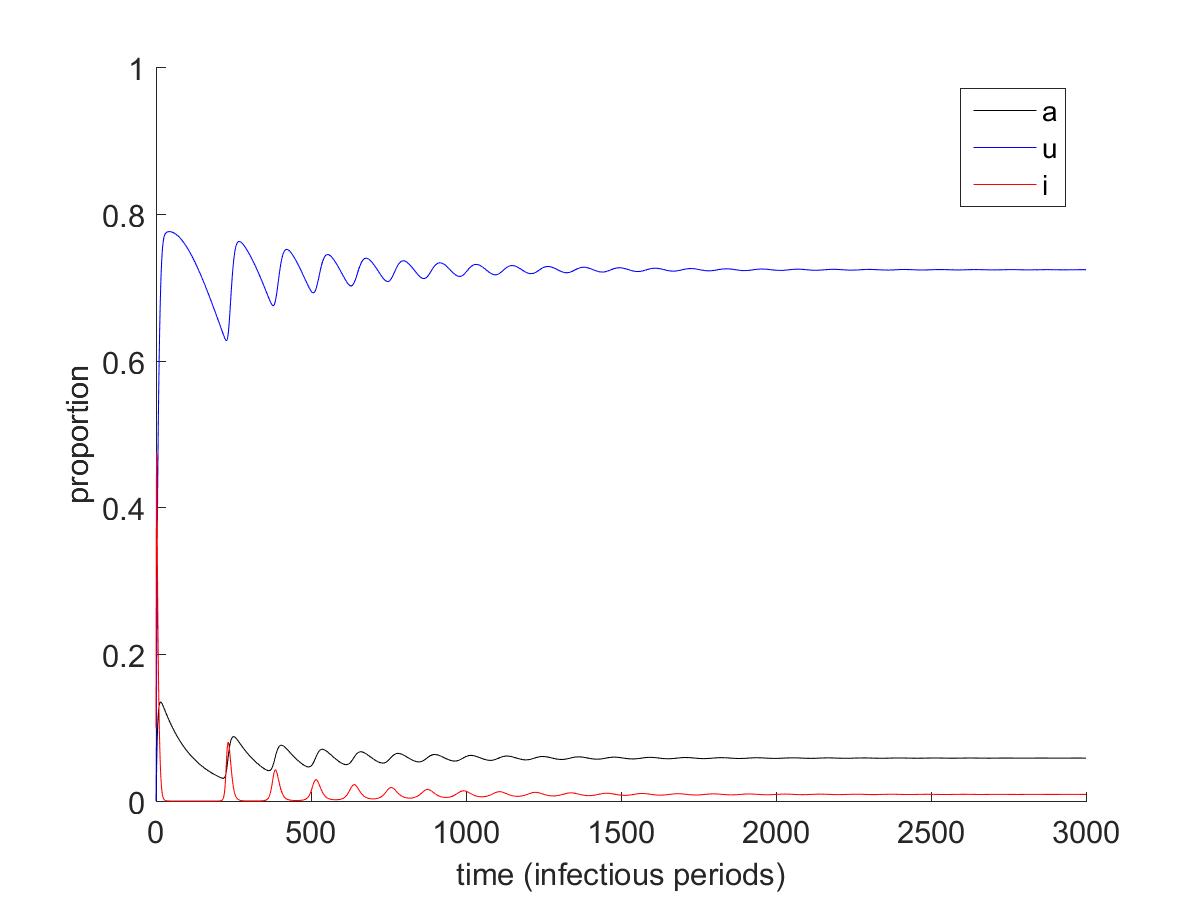}
&
\hspace{-0.75cm}

\includegraphics[scale = 0.2]{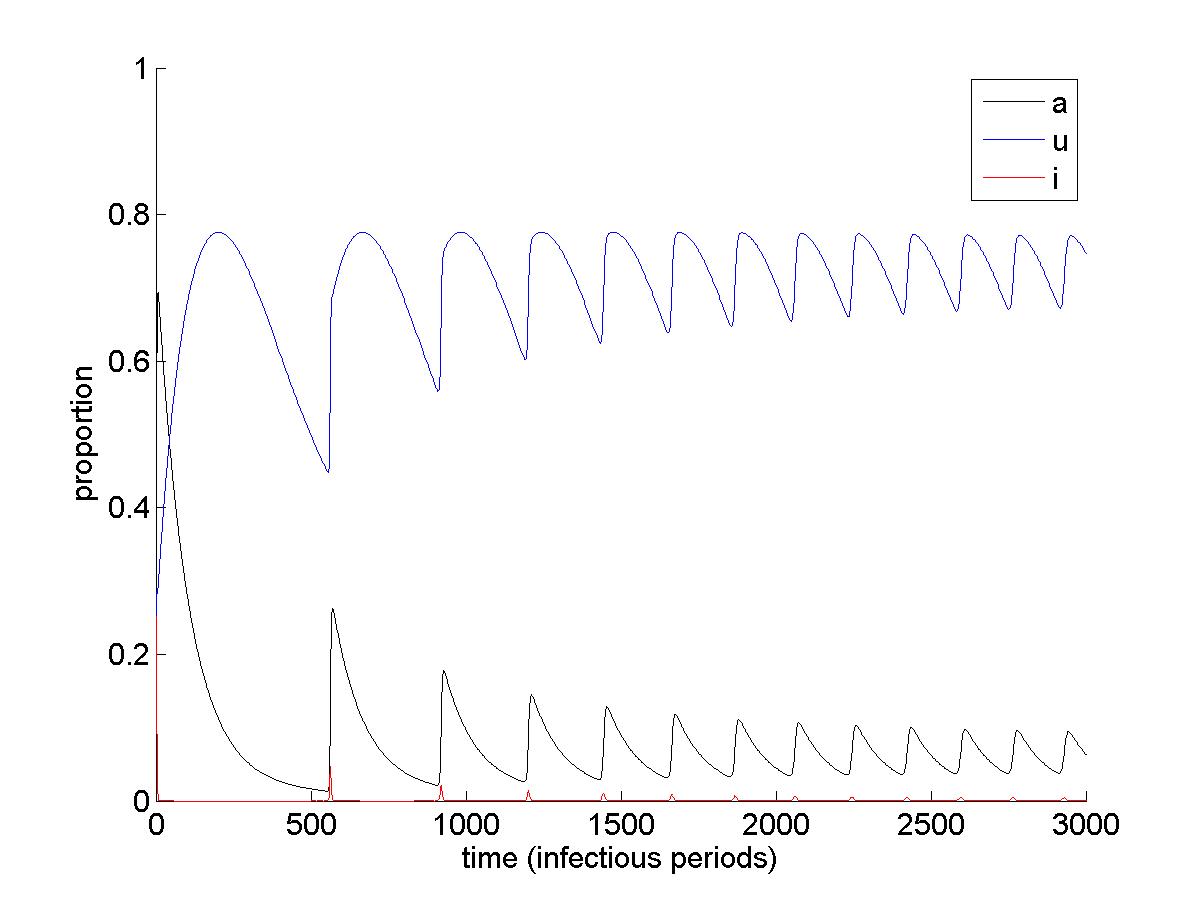}
\end{tabular}
\caption{Top two and bottom left: Evolution of the fraction of infectious, aware, and unwilling hosts for different values of $q$ along the vertical section in the Hopf-bifurcation curve corresponding to $p=0.05$ (when $p = 0.05$, $q^* = 0.1923$ is the only Hopf-bifurcation point): $q=0.05$ (top left), 0.15 (top right), 0.25 (bottom left). 
Bottom right: Evolution of the fraction of infectious, aware, and unwilling hosts for $p = 1$ and $q = 0$.
Initial condition:  $a(0) = u(0) = 0$, $i(0) = 0.1$.  Here, $\alpha_u = 3$, $\delta = 1$, $\delta_a = 0.01$, $\delta_u = 0.05$, $\beta = 3$, $\beta_u = 0.5$, $\alpha_a = 0.012$, $\beta_a=0.2$ and $\alpha_i=0.05$.
\label{Fig:Evolution-SAUIS}}
\end{center}
\end{figure}

\begin{figure}[H]
\begin{center}
\begin{tabular}{c}
\includegraphics[scale=0.2]{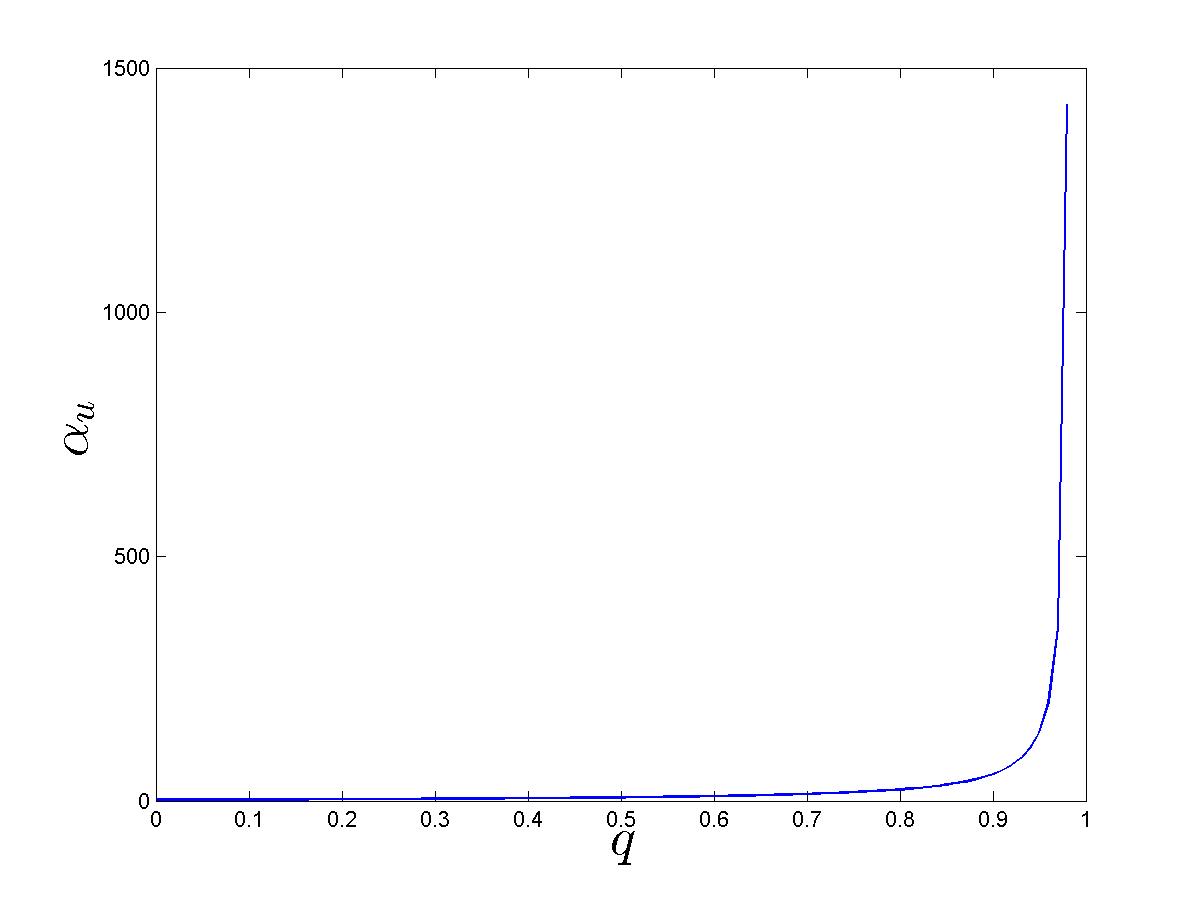}
\end{tabular}
\caption{Hopf-bifurcation curve. $q$: $0$(left) -- $0.98$(right). When $q = 0.1$, $\alpha_u^*=2.6449$ is the only Hopf-bifurcation point.  Here, $\delta = 1$, $\delta_a = 0.01$, $\delta_u = 0.05$, $\beta = 3$, $\beta_a = 0.2$, $\beta_u = 0.5$, $\alpha_a = 0.012$, $\alpha_i = 0.05$, $p = 1-q$.
}
\end{center}
\end{figure}

\begin{figure}[H]
\begin{center}
\begin{tabular}{cc}
\hspace{-1cm}
\includegraphics[scale=0.2]{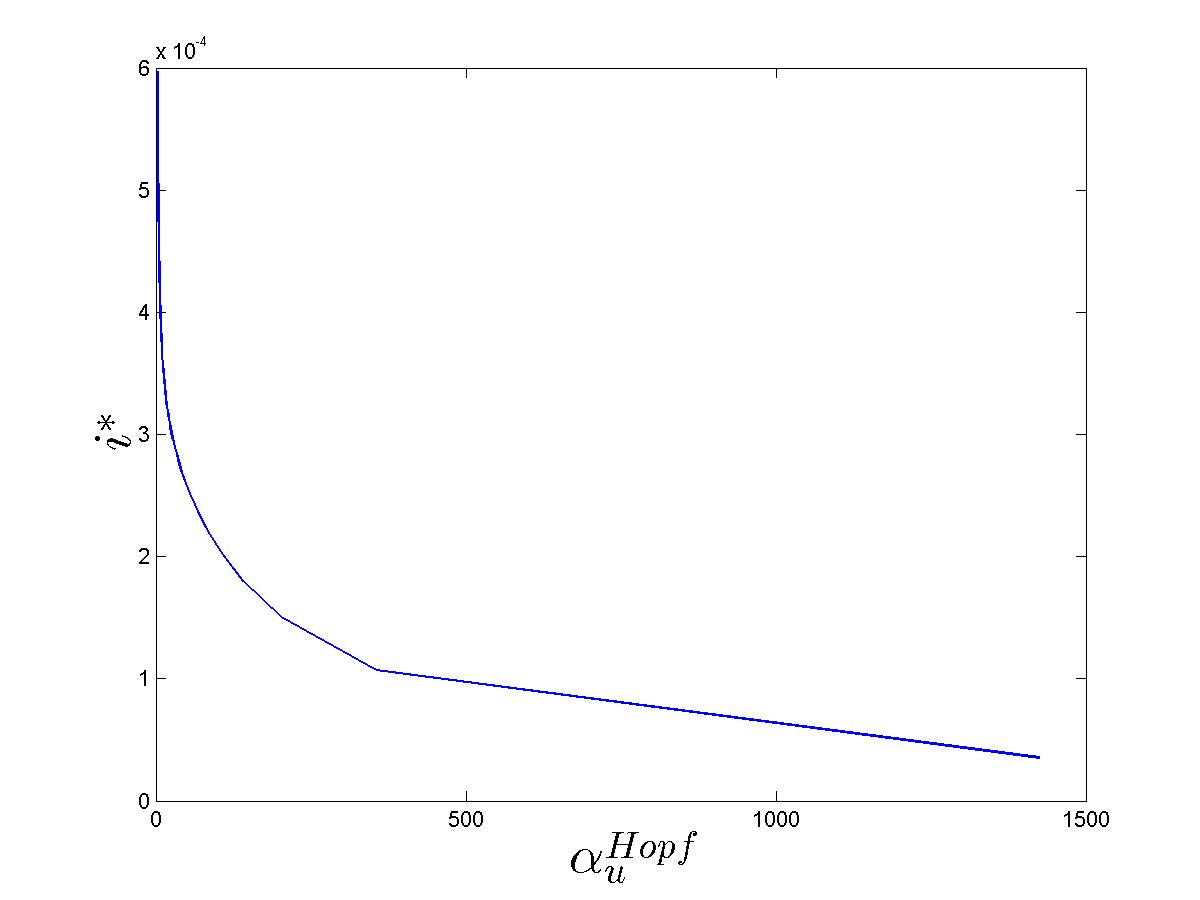}
\end{tabular}
\caption{
Fraction of infectious hosts as a function of $\alpha_u$ along the Hopf-bifurcation curve in Figure 11, where $\delta = 1$, $\delta_a = 0.01$, $\delta_u = 0.05$, $\beta = 3$, $\beta_a = 0.2$, $\beta_u = 0.5$, $\alpha_a = 0.012$, $\alpha_i = 0.05$, $p = 1-q$.
\label{Fig:Evolution-SAUIS}}
\end{center}
\end{figure}

\begin{figure}[H]
\begin{center}
\begin{tabular}{cc}
\hspace{-1cm}
\includegraphics[scale=0.2]{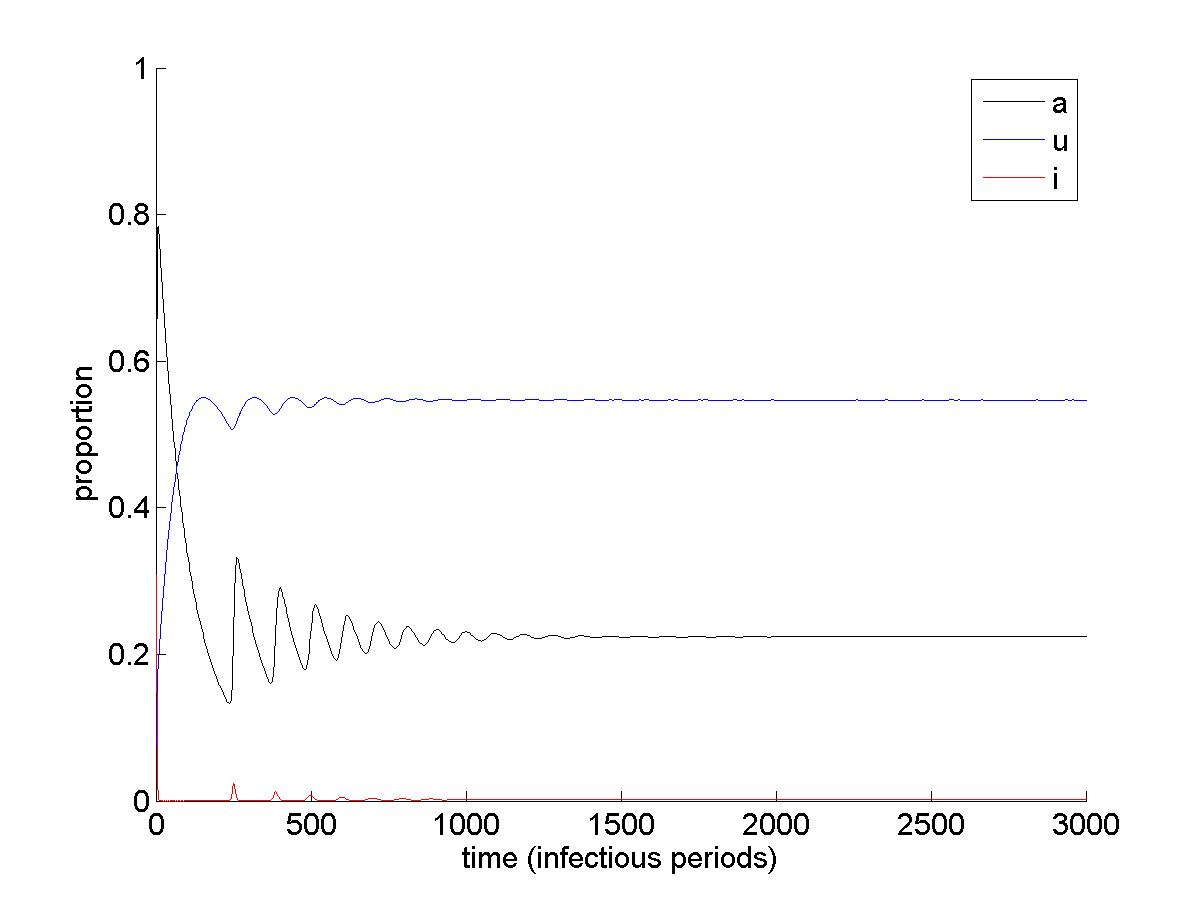}
&
\hspace{-0.75cm}
\includegraphics[scale=0.2]{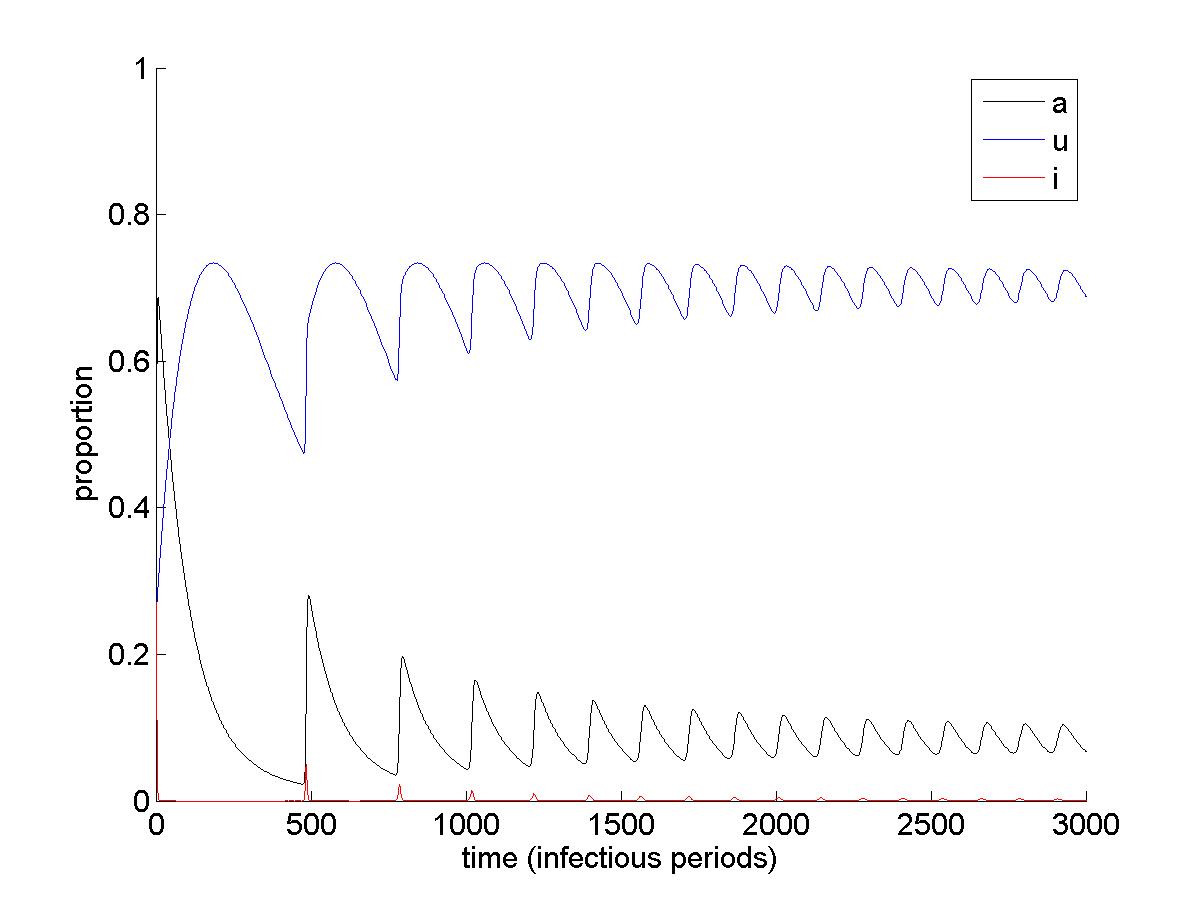}
\\
\hspace{-1cm}
\includegraphics[scale=0.2]{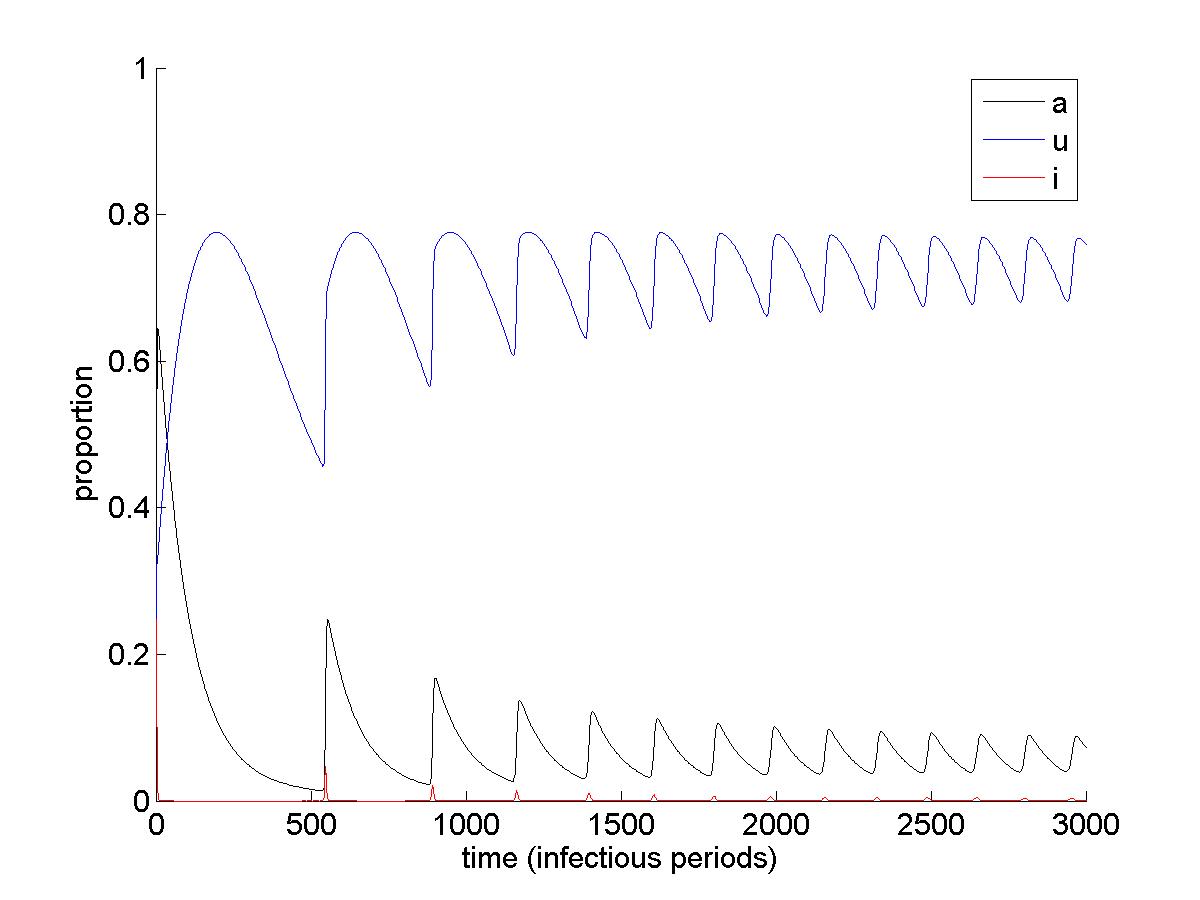}
&
\hspace{-0.75cm}
\includegraphics[scale=0.2]{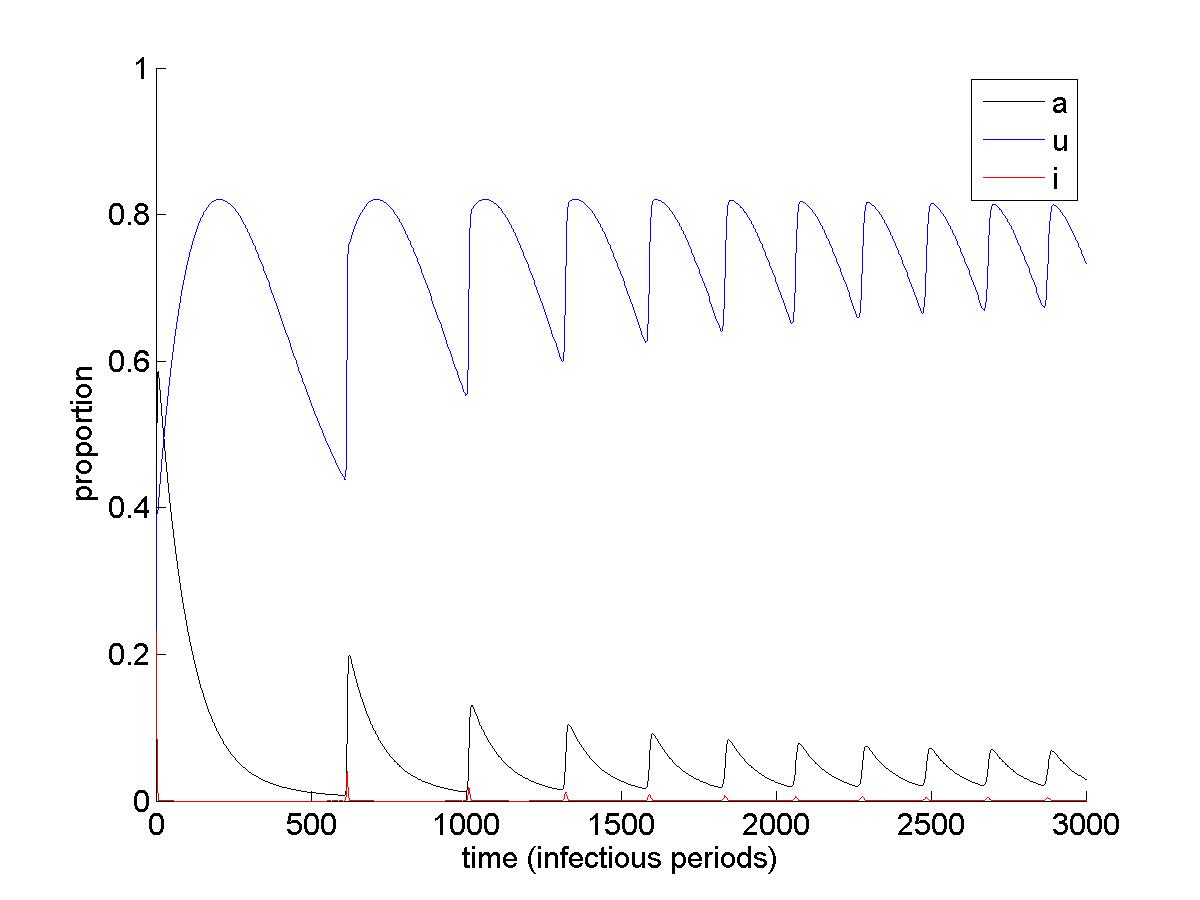}
\end{tabular}
\caption{Evolution of the fraction of infectious, aware, and unwilling hosts for different values of $\alpha_u$ along the vertical section in the Hopf-bifurcation curve corresponding to $q=0.1$ (when $q = 0.1$, $\alpha_u^* = 2.6449$ is the only Hopf-bifurcation point): $\alpha_u=0.5$ (top left), 2 (top right), 3 (bottom left), 5 (bottom right).  Initial condition:  $a(0) = u(0) = 0$, $i(0) = 0.1$.  Here $\delta = 1$, $\delta_a = 0.01$, $\delta_u = 0.05$, $\beta = 3$, $\beta_a = 0.2$, $\beta_u = 0.5$, $\alpha_a = 0.012$, $\alpha_i = 0.05$.
\label{Fig:Evolution-SAUIS}}
\end{center}
\end{figure}

\begin{figure}[H]
\begin{center}
\begin{tabular}{c}
\includegraphics[scale=0.2]{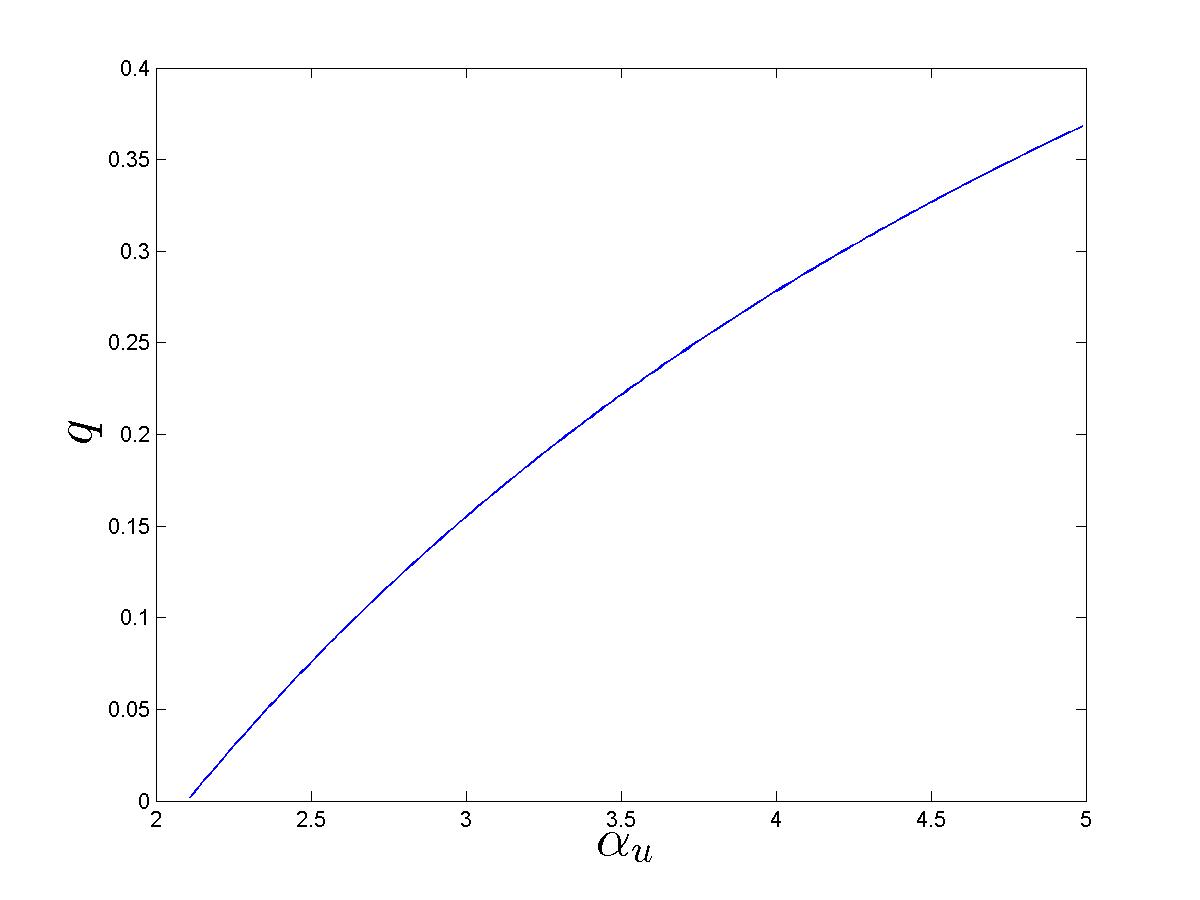}
\end{tabular}
\caption{Hopf-bifurcation curve. $\alpha_u$: $2.11$(left)-- right end. When $\alpha_u = 3$, $q^*=0.1549$ is the only Hopf-bifurcation point.  Here, $\delta = 1$, $\delta_a = 0.01$, $\delta_u = 0.05$, $\beta = 3$, $\beta_a = 0.2$, $\beta_u = 0.5$, $\alpha_a = 0.012$, $\alpha_i = 0.05$, $p = 1-q$.
}
\end{center}
\end{figure}

\begin{figure}[H]
\begin{center}
\begin{tabular}{cc}
\hspace{-1cm}
\includegraphics[scale=0.2]{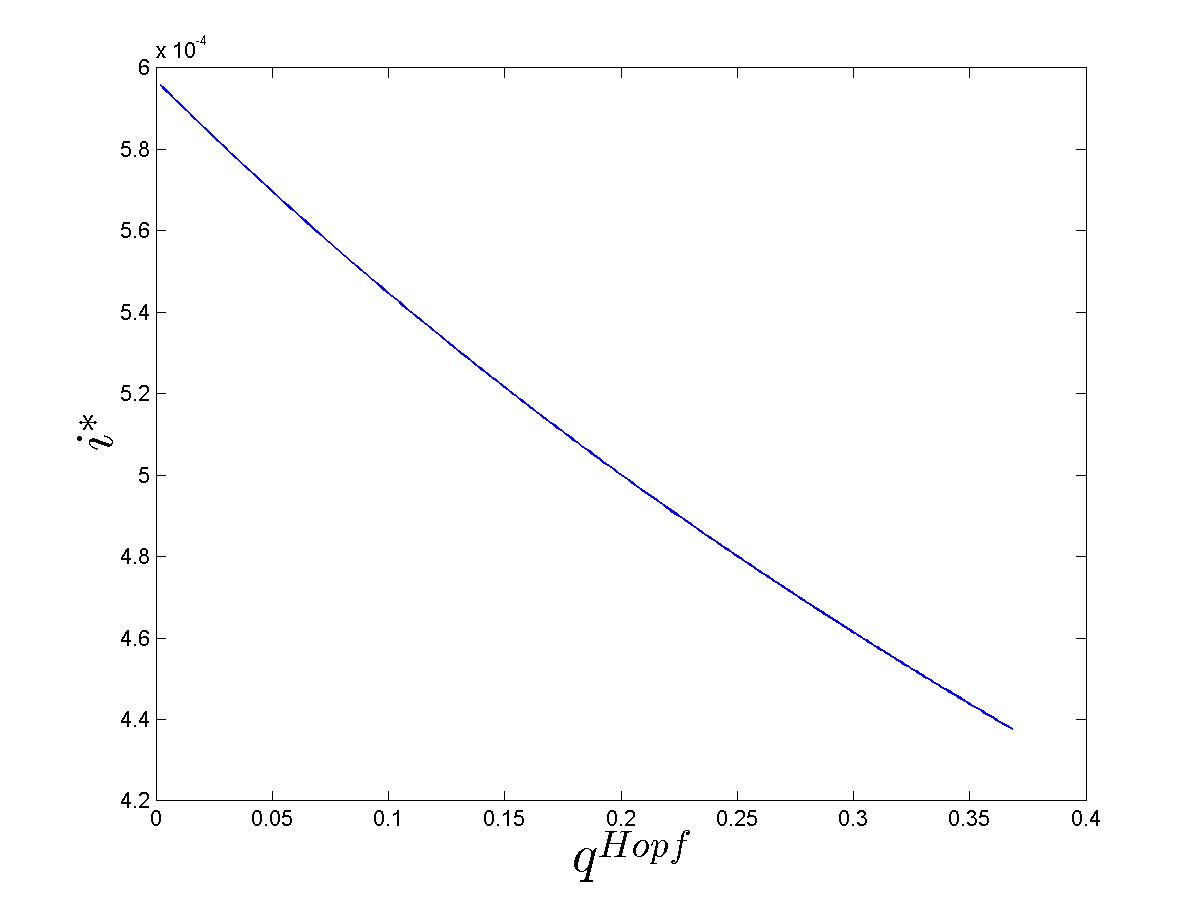}
\end{tabular}
\caption{
Fraction of infectious hosts as a function of $q$ along the Hopf-bifurcation curve in Figure 14, where $\delta = 1$, $\delta_a = 0.01$, $\delta_u = 0.05$, $\beta = 3$, $\beta_a = 0.2$, $\beta_u = 0.5$, $\alpha_a = 0.012$, $\alpha_i = 0.05$, $p = 1-q$.
}
\end{center}
\end{figure}

\begin{figure}[H]
\begin{center}
\begin{tabular}{cc}
\hspace{-1cm}
\includegraphics[scale=0.2]{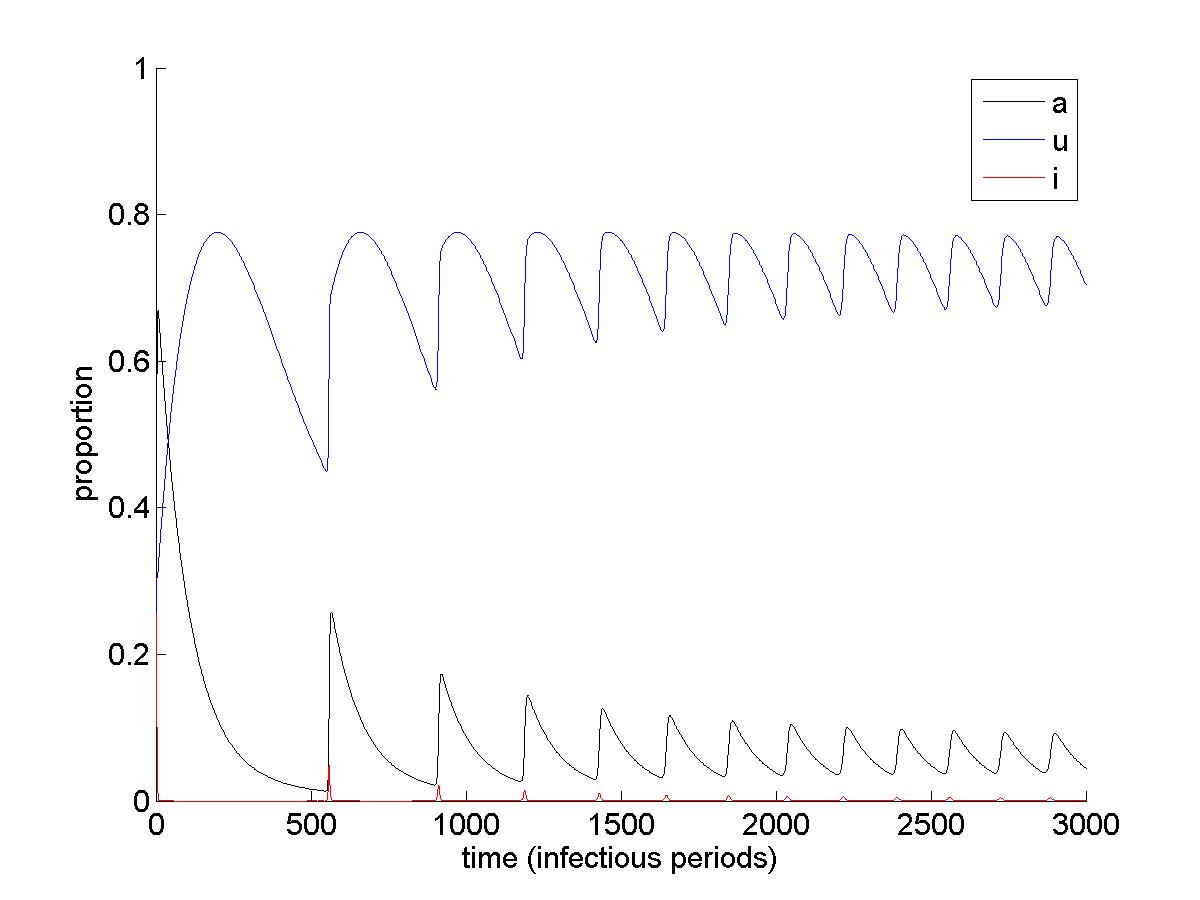}
&
\hspace{-0.75cm}
\includegraphics[scale=0.2]{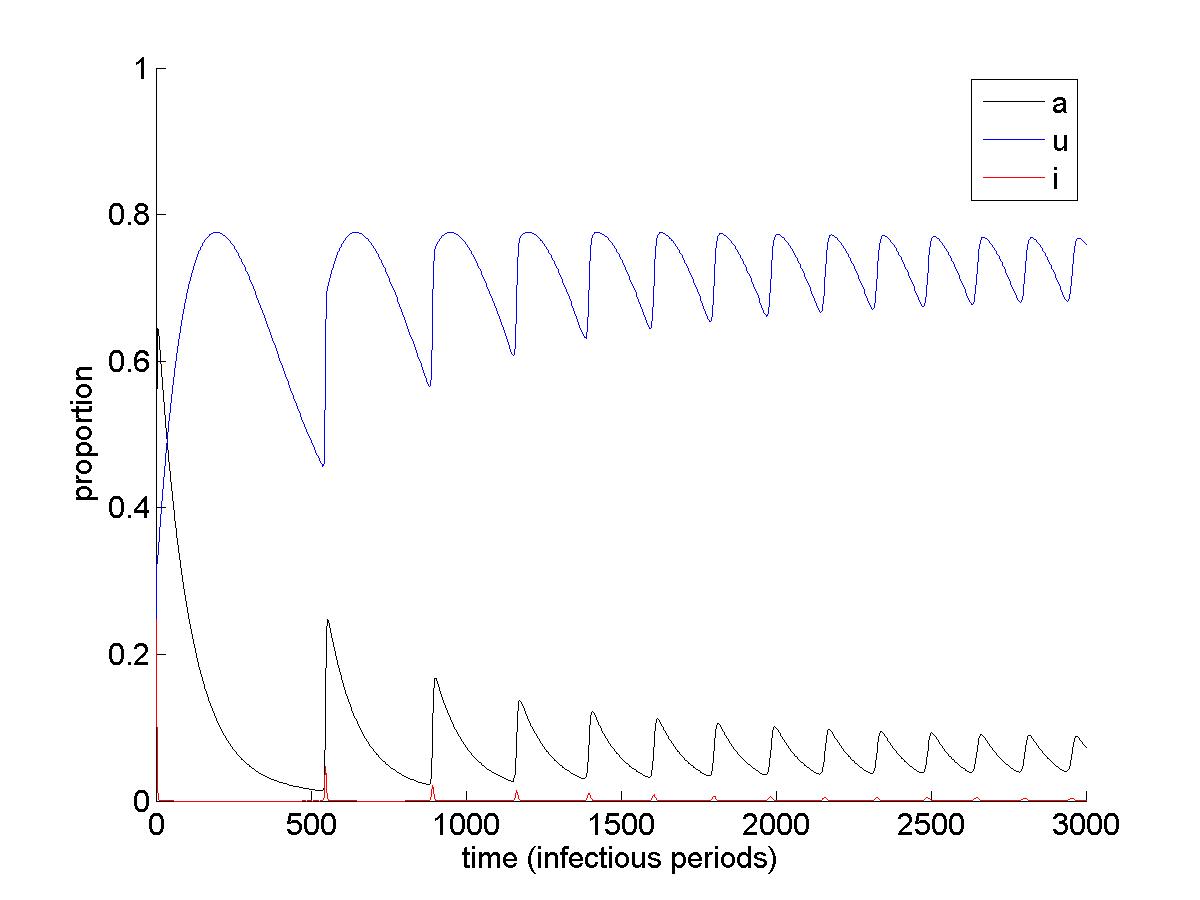}
\\
\hspace{-1cm}
\includegraphics[scale=0.2]{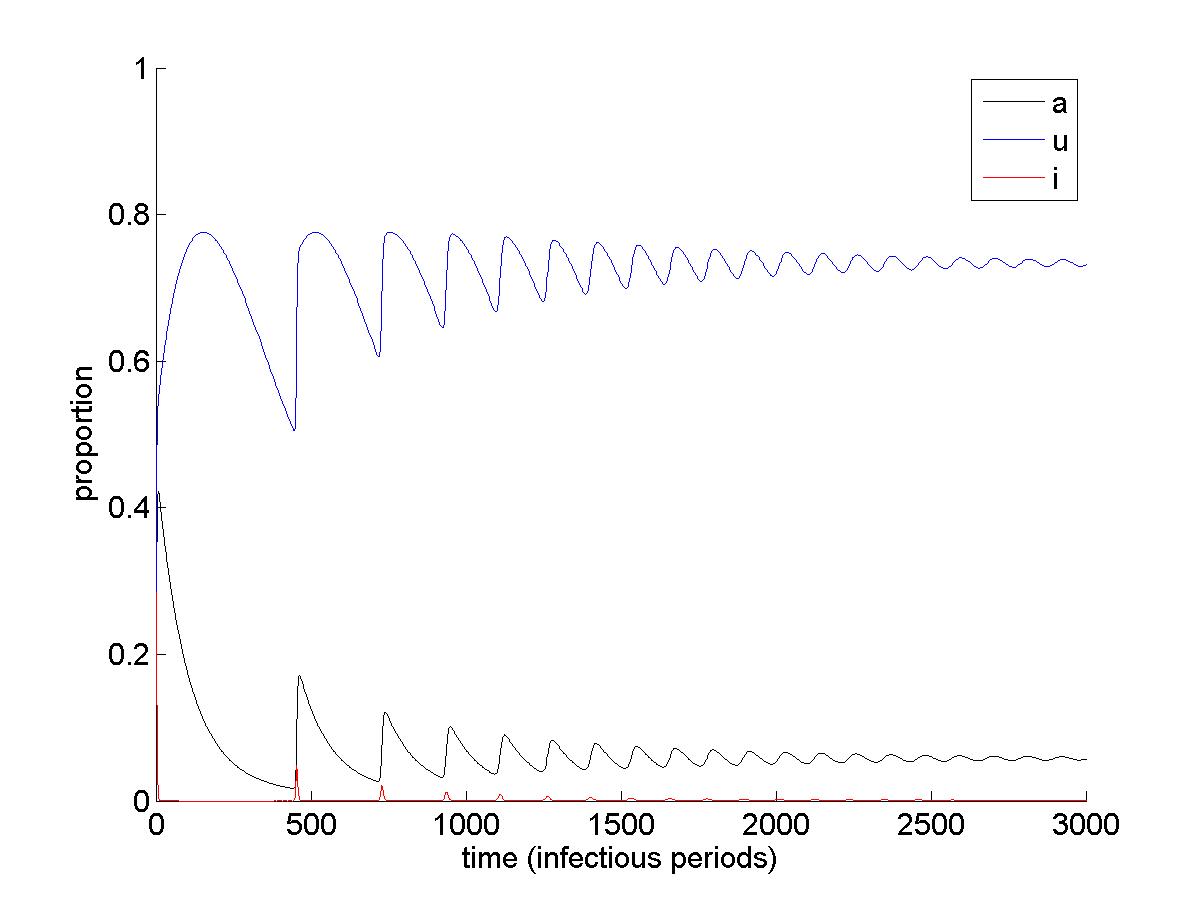}
&
\hspace{-0.75cm}
\includegraphics[scale=0.2]{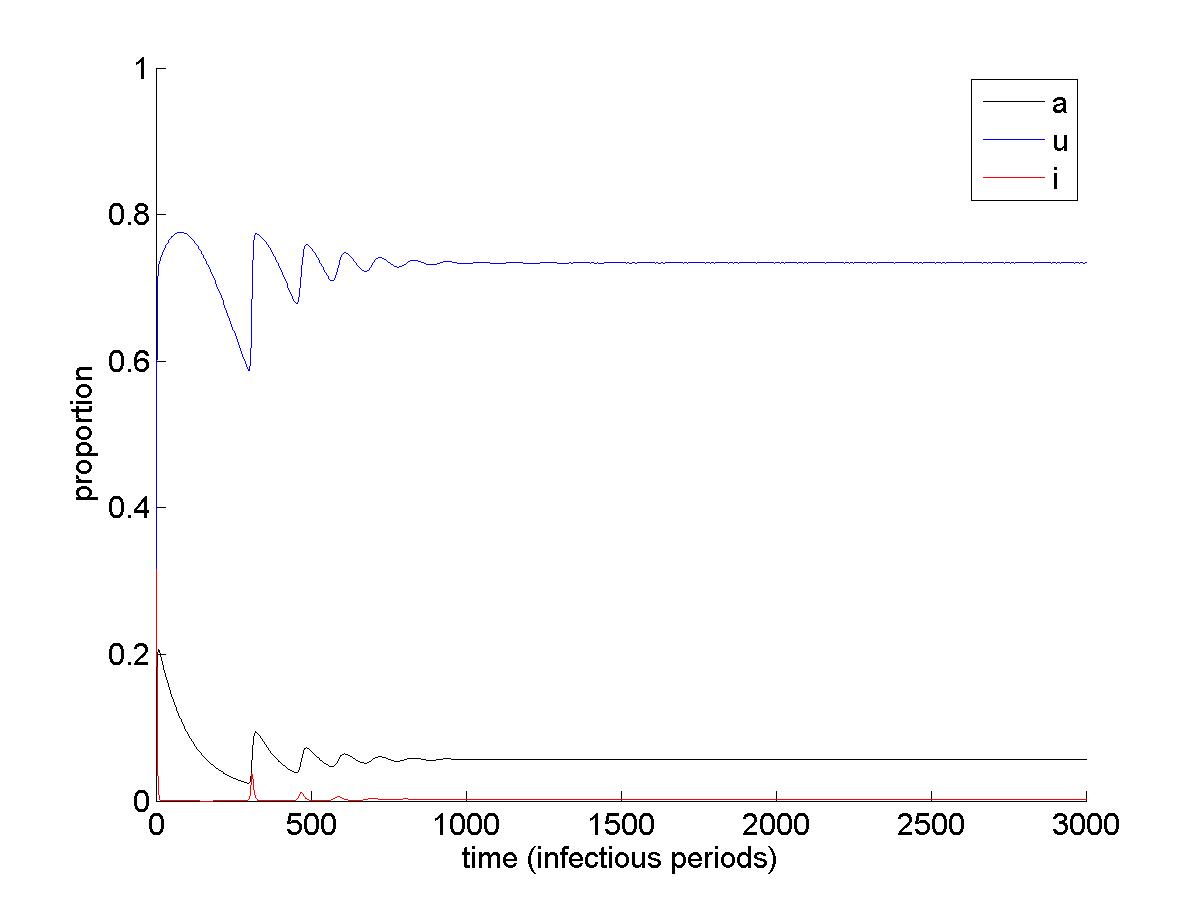}
\end{tabular}
\caption{Evolution of the fraction of infectious, aware, and unwilling hosts for different values of $q$ along the vertical section in the Hopf-bifurcation curve corresponding to $\alpha_u=3$ (when $\alpha_u = 3$, $q^* = 0.1549$ is the only Hopf-bifurcation point): $q=0.05$ (top left), 0.1 (top right), 0.5 (bottom left), 0.8 (bottom right).  Initial condition:  $a(0) = u(0) = 0$, $i(0) = 0.1$.  Here, $\delta = 1$, $\delta_a = 0.01$, $\delta_u = 0.05$, $\beta = 3$, $\beta_a = 0.2$, $\beta_u = 0.5$, $\alpha_a = 0.012$, $\alpha_i = 0.05$.
\label{Fig:Evolution-SAUIS}}
\end{center}
\end{figure}

\begin{figure}[H]
\begin{center}
\begin{tabular}{c}
\includegraphics[scale=0.2]{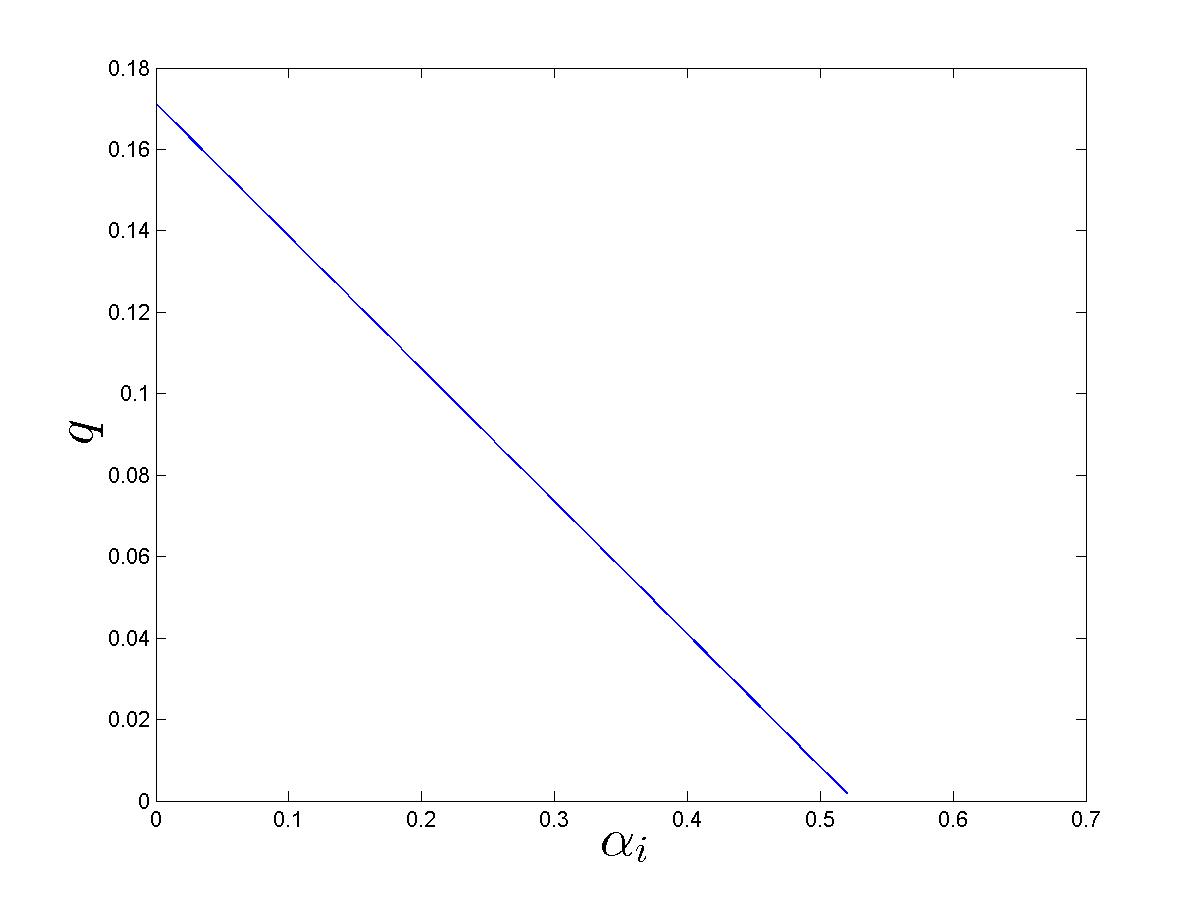}
\end{tabular}
\caption{Hopf-bifurcation curve. $\alpha_i$: $0$(left) -- $0.52$(right). When $\alpha_i = 0.1$, $q^*=0.1386$ is the only Hopf-bifurcation point.  Here, $\delta = 1$, $\delta_a = 0.01$, $\delta_u = 0.05$, $\beta = 3$, $\beta_a = 0.2$, $\beta_u = 0.5$, $\alpha_a = 0.012$, $\alpha_u = 1$, $p = 1-q$.
}
\end{center}
\end{figure}

\begin{figure}[H]
\begin{center}
\begin{tabular}{cc}
\hspace{-1cm}
\includegraphics[scale=0.2]{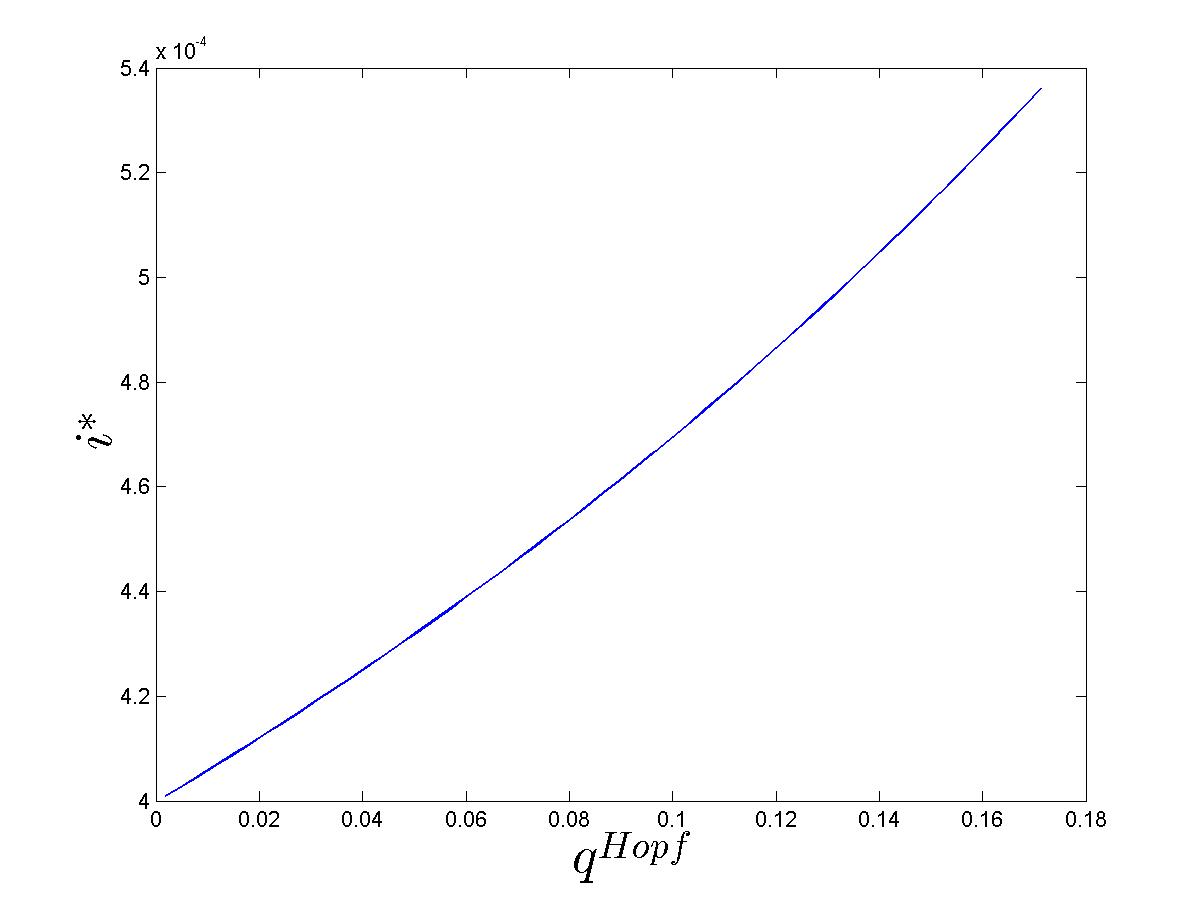}
\end{tabular}
\caption{ 
Fraction of infectious hosts as a function of $q$ along the Hopf-bifurcation curve in Figure 17, where $\delta = 1$, $\delta_a = 0.01$, $\delta_u = 0.05$, $\beta = 3$, $\beta_a = 0.2$, $\beta_u = 0.5$, $\alpha_a = 0.012$, $\alpha_u = 1$, $p = 1-q$.
}
\end{center}
\end{figure}

\begin{figure}[H]
\begin{center}
\begin{tabular}{cc}
\hspace{-1cm}
\includegraphics[scale=0.2]{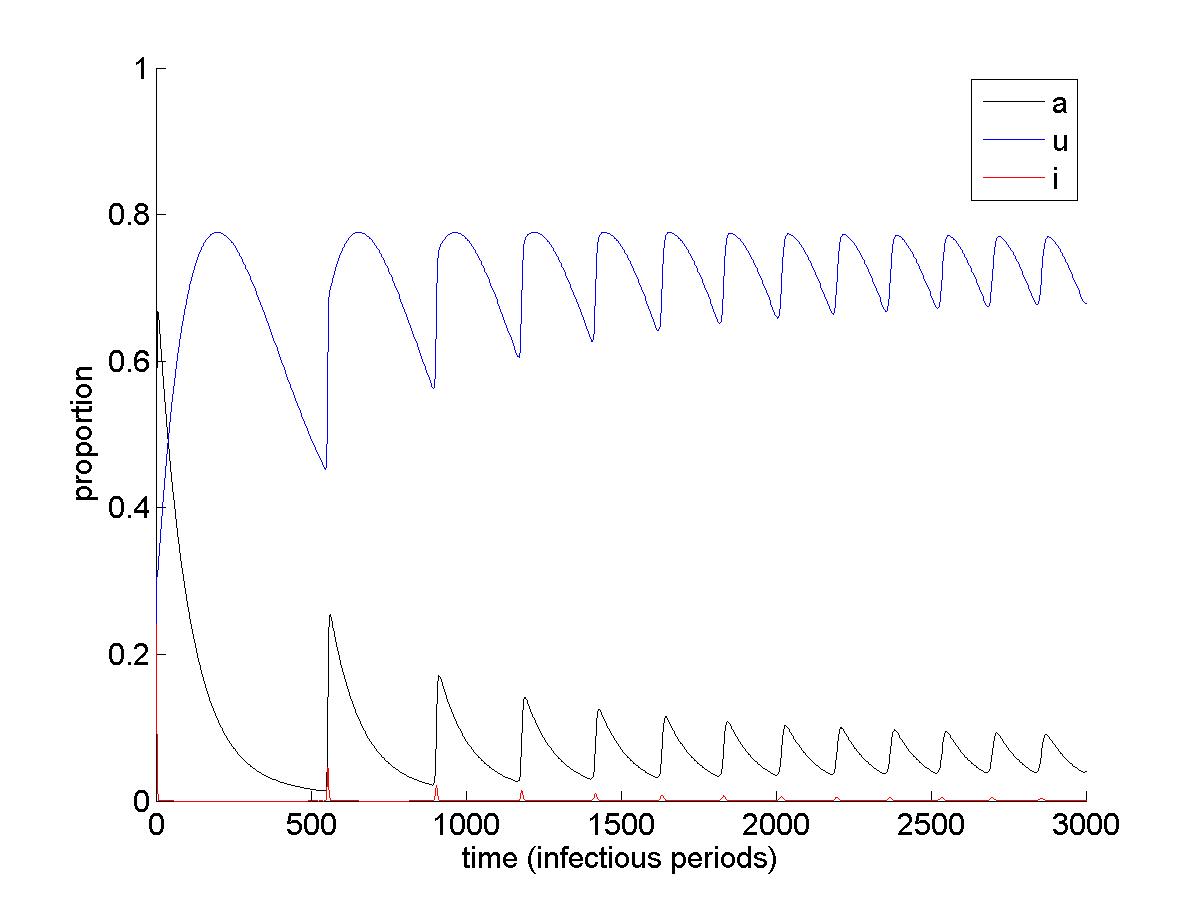}
&
\hspace{-0.75cm}
\includegraphics[scale=0.2]{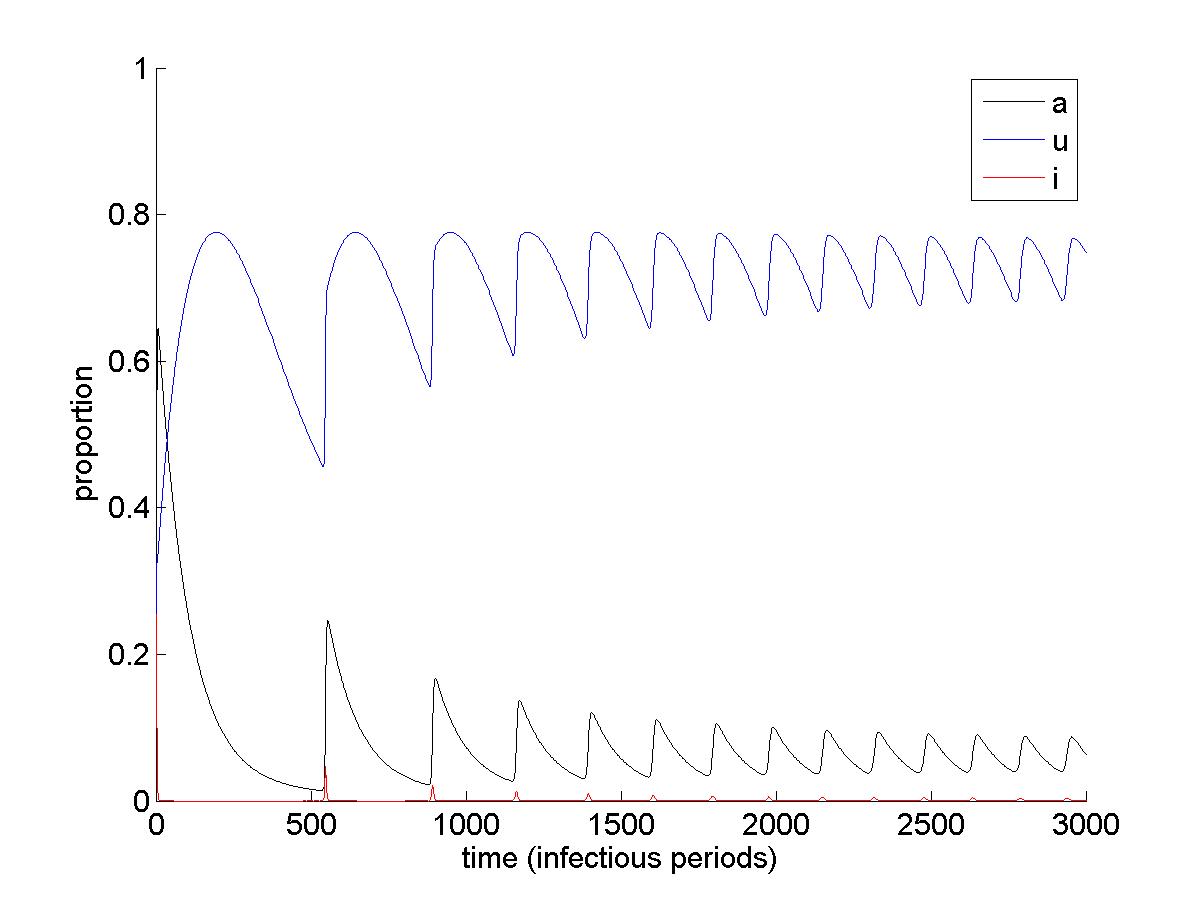}
\\
\hspace{-1cm}
\includegraphics[scale=0.2]{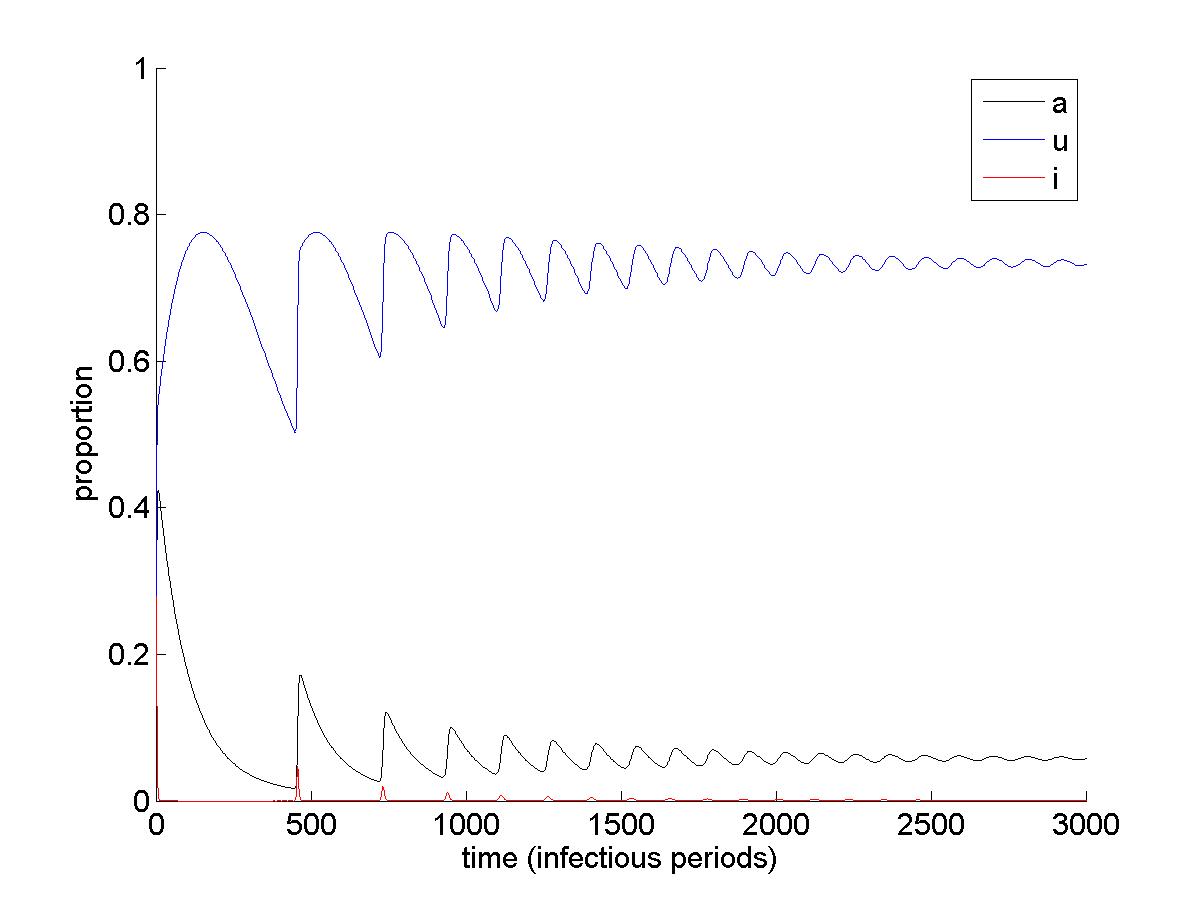}
&
\hspace{-0.75cm}
\includegraphics[scale=0.2]{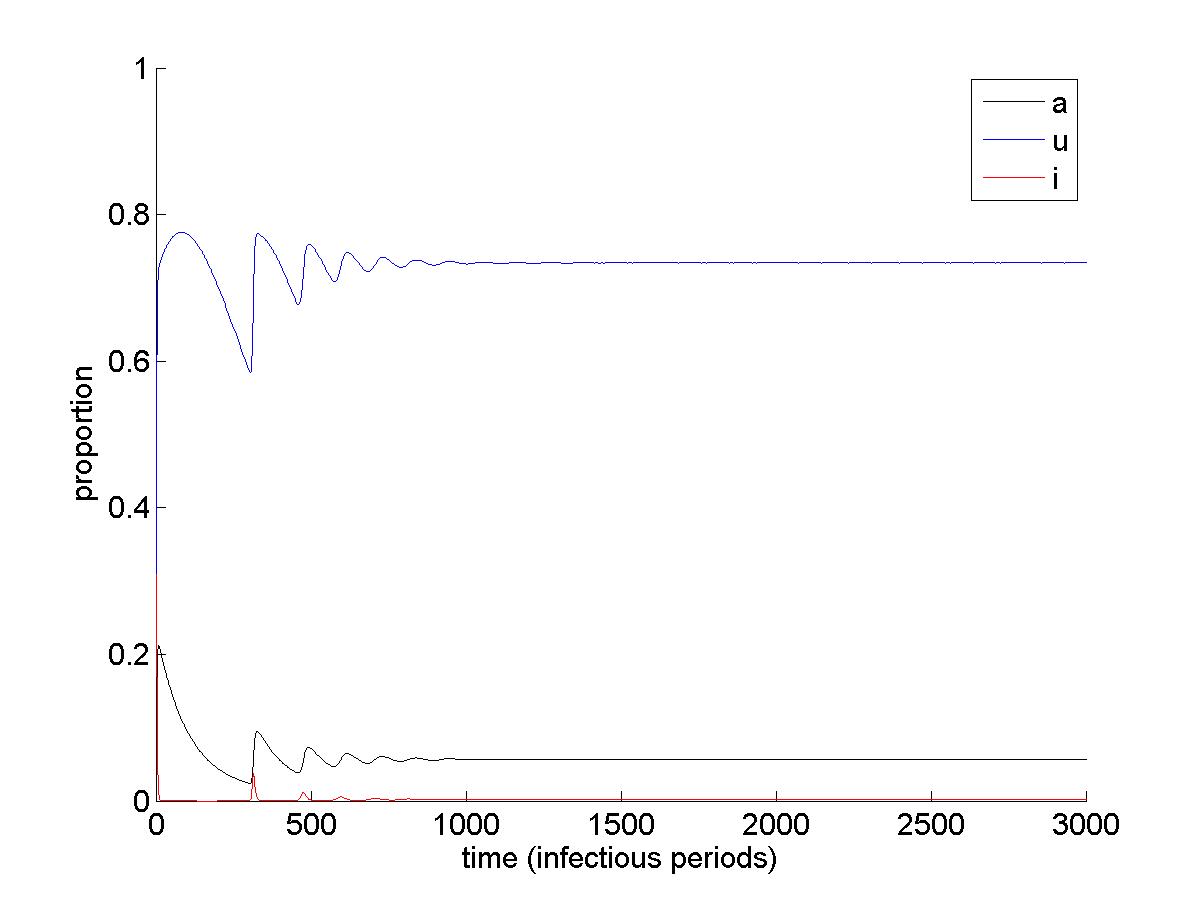}
\end{tabular}
\caption{Evolution of the fraction of infectious, aware, and unwilling hosts for different values of $q$ along the vertical section in the Hopf-bifurcation curve corresponding to $\alpha_i=0.1$ (when $\alpha_i = 0.1$, $q^* = 0.1386$ is the only Hopf-bifurcation point): $q=0.05$ (top left), 0.1 (top right), 0.5 (bottom left), 0.8 (bottom right).  Initial condition:  $a(0) = u(0) = 0$, $i(0) = 0.1$.  Here, $\delta = 1$, $\delta_a = 0.01$, $\delta_u = 0.05$, $\beta = 3$, $\beta_a = 0.2$, $\beta_u = 0.5$, $\alpha_a = 0.012$, and $\alpha_u = 3$.
\label{Fig:Evolution-SAUIS}}
\end{center}
\end{figure}

\begin{figure}[H]
\begin{center}
\begin{tabular}{c}
\includegraphics[scale=0.2]{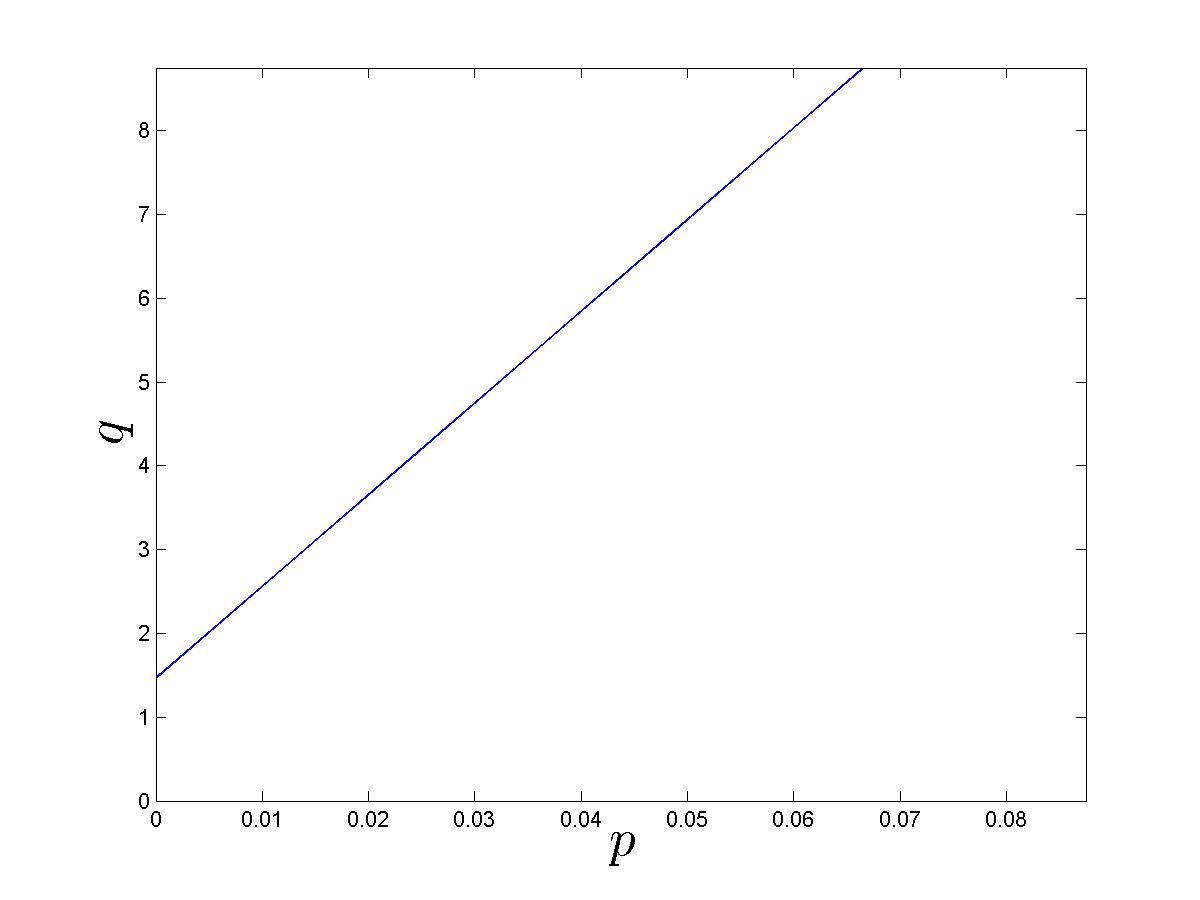}
\end{tabular}
\caption{Hopf-bifurcation curve where $(p, q)$ is the Hopf pair.  Only $(p, q)$ pairs with $p + q \leq 1$ make biological sense.  As sustained oscillations can occur in the region below the curve, the figure show that this will be the case in the entire biologically feasible region with parameter settings $\beta = 3$, $\beta_a = 0.2$, $\beta_u = 0.4$, $\alpha_i = 0.05$, $\alpha_a = 0.012$, $\delta_a = 0.01$, $\delta_u = 0.05$, $\delta = 1.7$, and $\alpha_u = 30$.
}
\end{center}
\end{figure}

\begin{figure}[H]
\begin{center}
\begin{tabular}{c}
\includegraphics[scale=0.2]{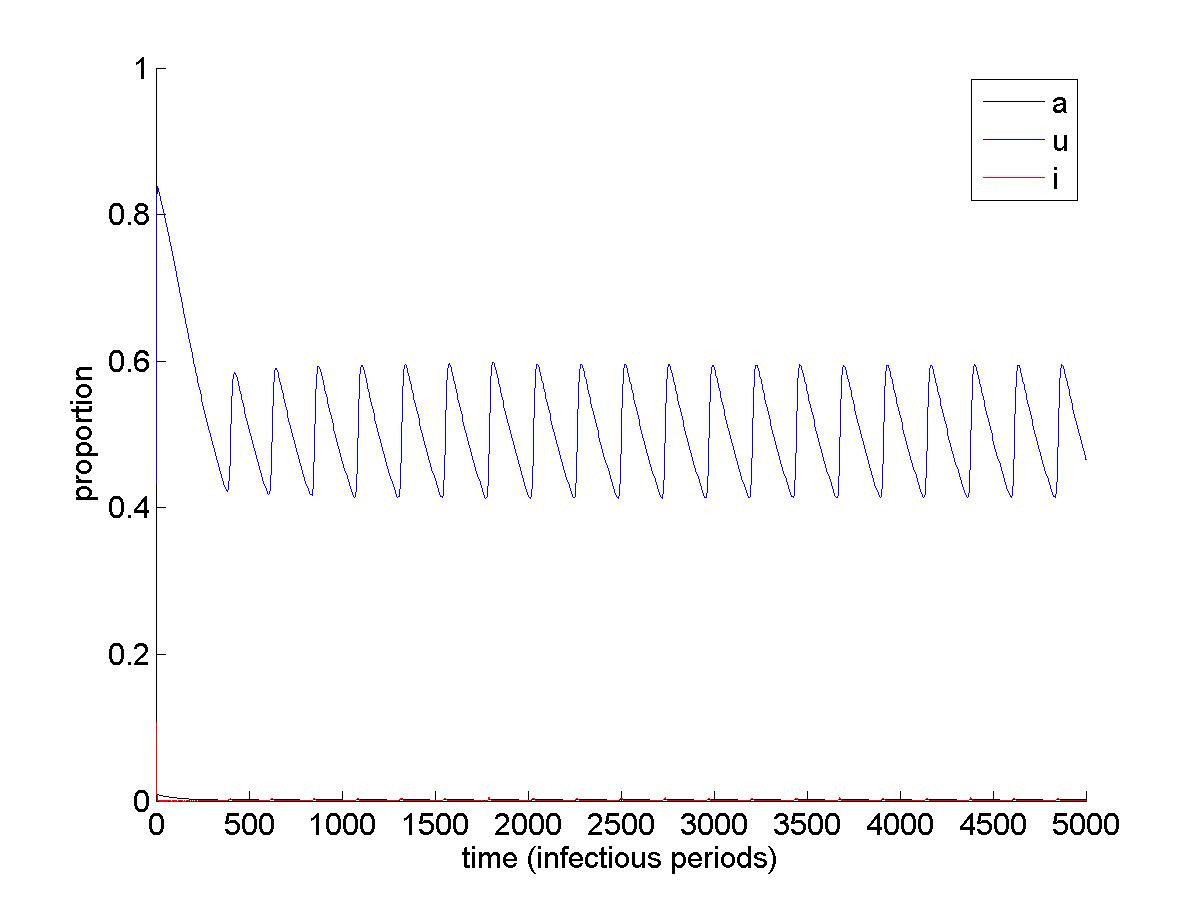}
\end{tabular}
\caption{Time series with $p = 0$, $q = 1$, $\beta = 3$, $\beta_a = 0.2$, $\beta_u = 0.4$, $\alpha_i = 0.05$, $\alpha_a = 0.012$, $\delta_a = 0.01$, $\delta_u = 0.05$, $\delta = 1.7$, and $\alpha_u = 30$.
}
\end{center}
\end{figure}

\begin{figure}[H]
\begin{center}
\begin{tabular}{c}
\includegraphics[scale=0.2]{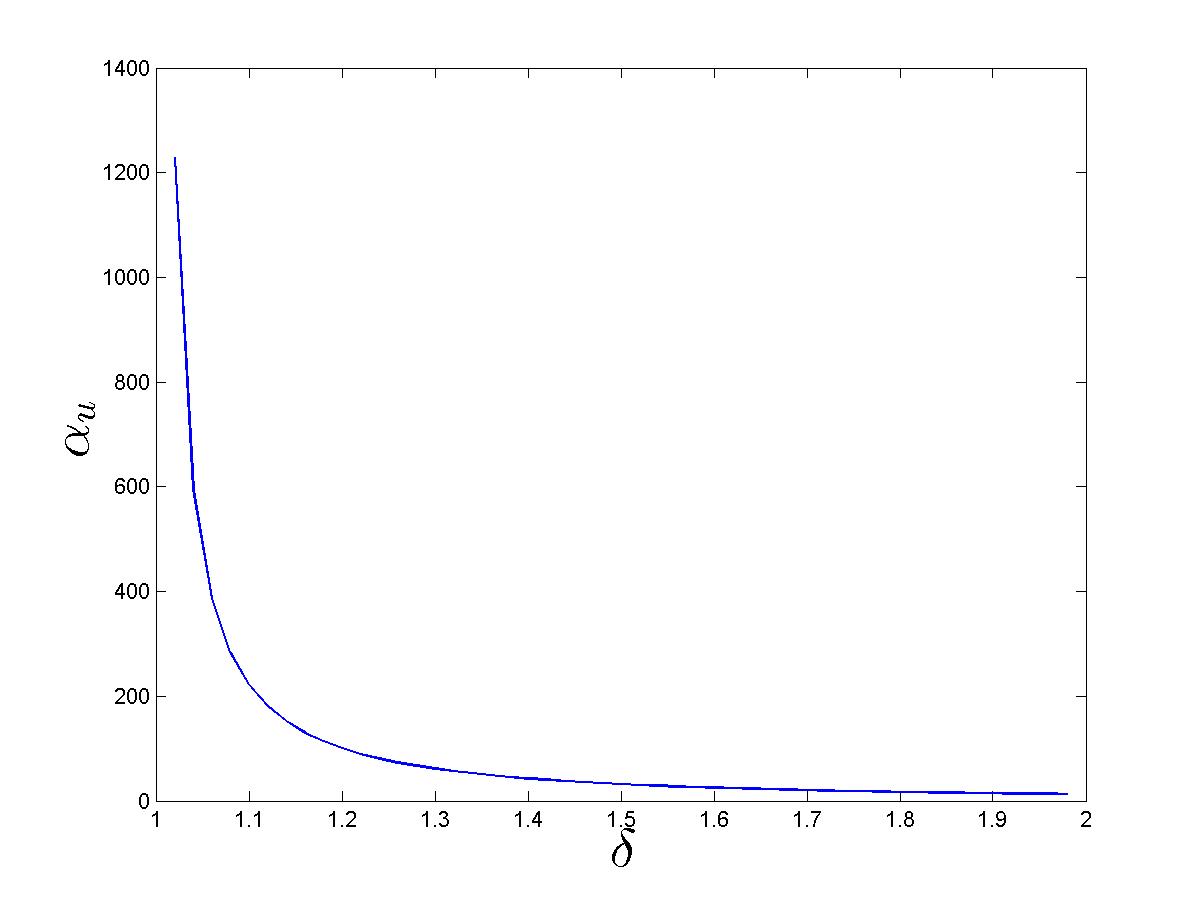}
\end{tabular}
\caption{Hopf-bifurcation curve where $(\delta, \alpha_u)$ is taken as the Hopf pair.  Here, $p = 0$, $q = 1$, $\beta = 3$, $\beta_a = 0.2$, $\beta_u = 0.4$, $\alpha_i = 0.05$, $\alpha_a = 0.012$, $\delta_a = 0.01$, and $\delta_u = 0.05$. The left end of $\delta$ is $1.02$.  When $\delta = 1.7$, $\alpha_u = 20.36$ is the only Hopf-bifurcation point.
}
\end{center}
\end{figure}

\begin{figure}[H]
\begin{center}
\begin{tabular}{cc}
\hspace{-1cm}
\includegraphics[scale=0.2]{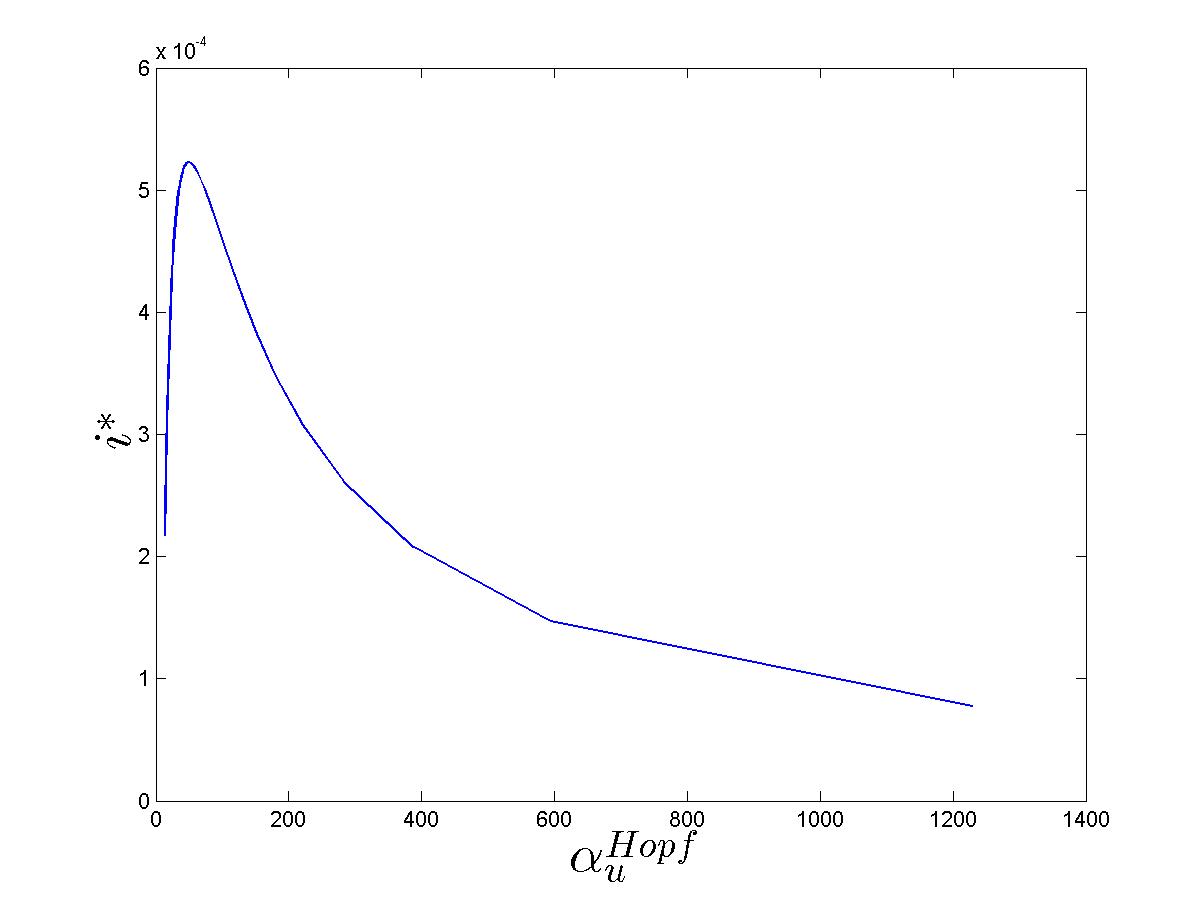}
\end{tabular}
\caption{
Fraction of infectious hosts as a function of $\alpha_u$ along the Hopf-bifurcation curve in Figure 22, where $p = 0$, $q = 1$, $\beta = 3$, $\beta_a = 0.2$, $\beta_u = 0.4$, $\alpha_i = 0.05$, $\alpha_a = 0.012$, $\delta_a = 0.01$, and $\delta_u = 0.05$.
}
\end{center}
\end{figure}

\begin{figure}[H]
\begin{center}
\begin{tabular}{c}
\includegraphics[scale=0.2]{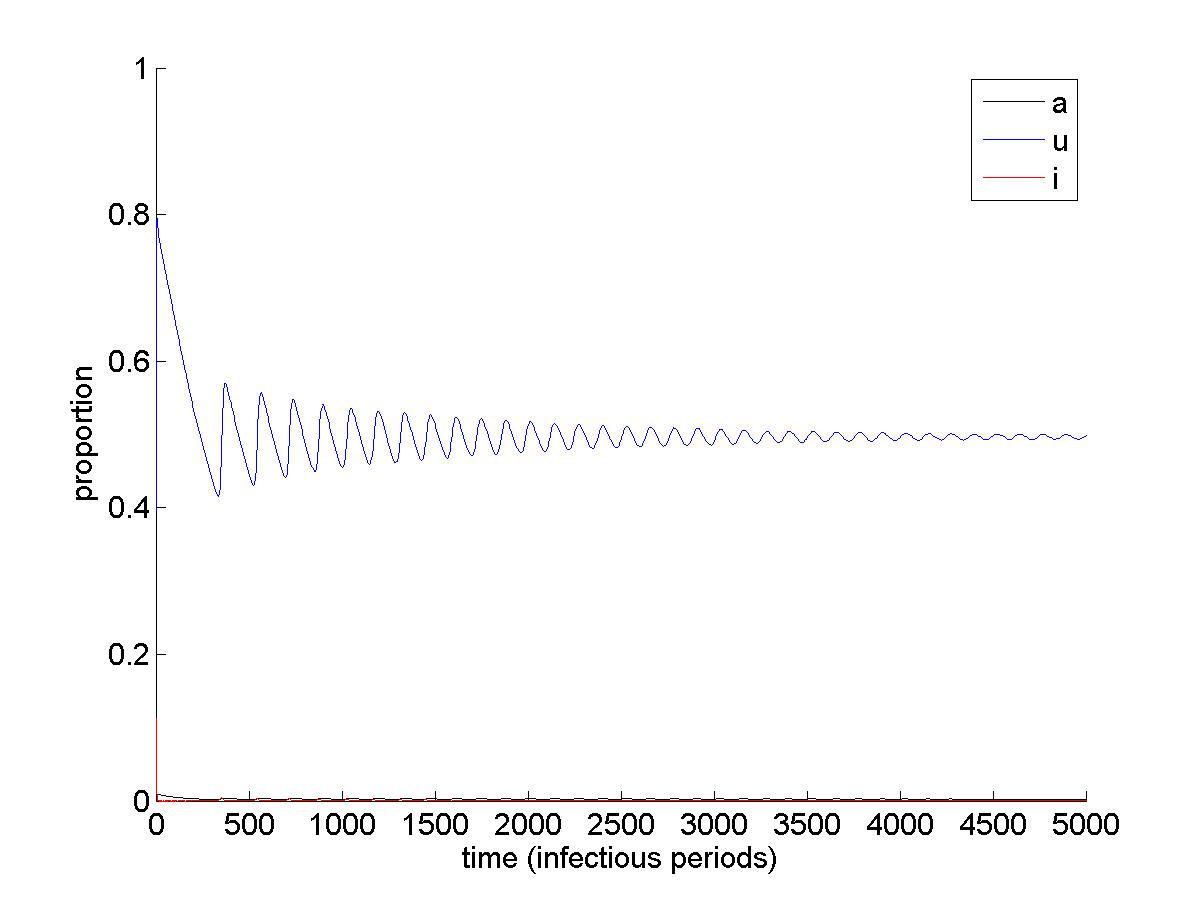}
\end{tabular}
\caption{Time series with $p = 0$, $q = 1$, $\beta = 3$, $\beta_a = 0.2$, $\beta_u = 0.4$, $\alpha_i = 0.05$, $\alpha_a = 0.012$, $\delta_a = 0.01$, $\delta_u = 0.05$, $\delta = 1.7$ and $\alpha_u = 19$.
}
\end{center}
\end{figure}

\begin{figure}[H]
\begin{center}
\begin{tabular}{c}
\includegraphics[scale=0.2]{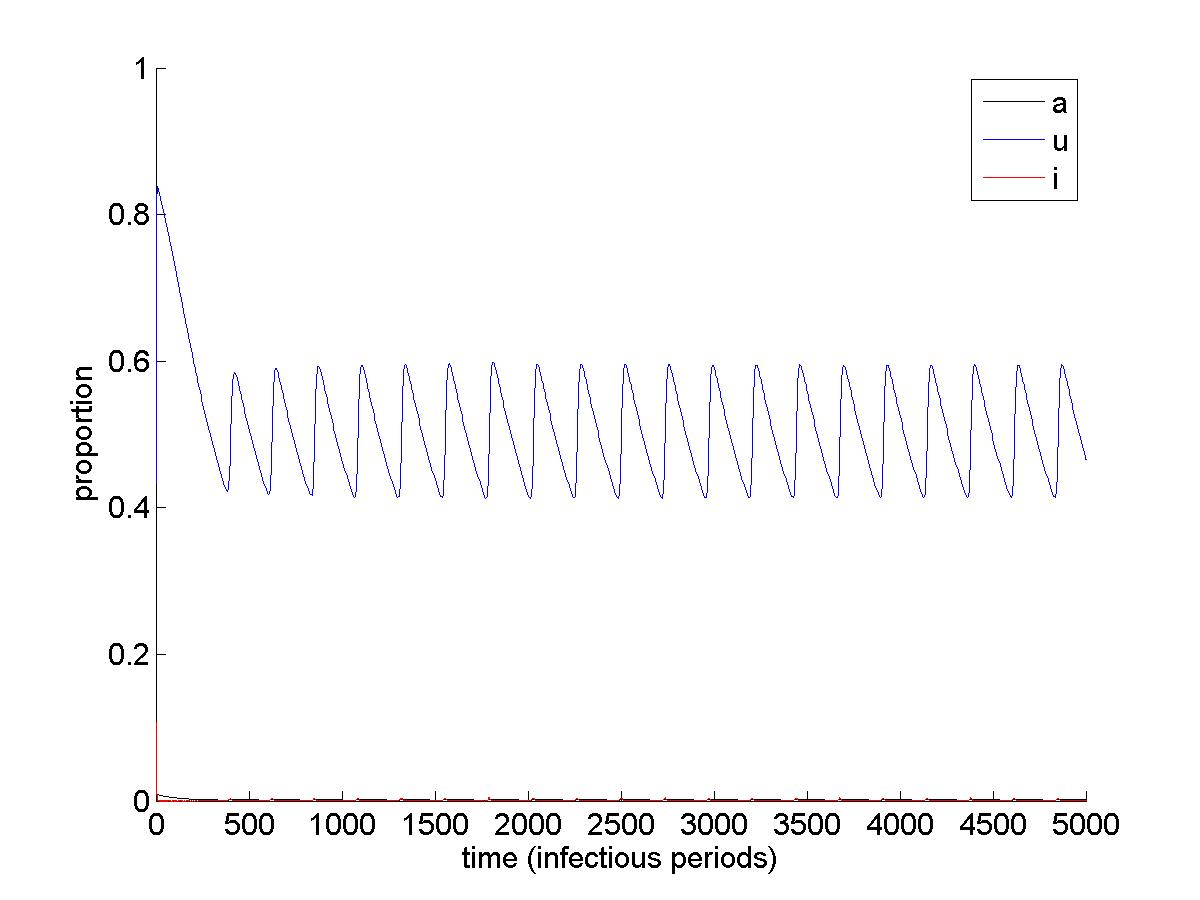}
\end{tabular}
\caption{Time series with $p = 0$, $q = 1$, $\beta = 3$, $\beta_a = 0.2$, $\beta_u = 0.4$, $\alpha_i = 0.05$, $\alpha_a = 0.012$, $\delta_a = 0.01$, $\delta_u = 0.05$, $\delta = 1.7$ and $\alpha_u = 30$.
}
\end{center}
\end{figure}

\end{document}